\documentclass[11pt]{amsart}
\usepackage{amsfonts}
\usepackage{amsmath}
\usepackage{amsthm}
\usepackage{bm}
\usepackage[dvips]{graphicx}
\usepackage{caption}
\usepackage{color}

\usepackage{bigints}

\numberwithin{equation}{section}
\allowdisplaybreaks

\theoremstyle{plain}
\newtheorem{theorem}{Theorem}[section]
\newtheorem{lemma}[theorem]{Lemma}
\newtheorem{proposition}[theorem]{Proposition}

\theoremstyle{definition}
\newtheorem{definition}[theorem]{Definition}
\newtheorem{remark}[theorem]{Remark}
\newtheorem{example}[theorem]{Example}

\def\beqn{\begin{equation}}
\def\beqn*{$$}
\def\eeqn{\end{equation}}
\def\eeqn*{$$}

\newcommand{\BX}{{\bf X}}

\newcommand{\BV}{{\bf V}}
\newcommand{\BW}{{\bf W}}
\newcommand{\bx}{{\bf x}}
\newcommand{\by}{{\bf y}}
\newcommand{\bz}{{\bf z}}

\def\P{\mathbb{P}}
\def\E{\mathbb{E}}
\def\Pn{\mathcal P_n}

\newcommand{\reals}{{\mathbb R}}
\newcommand{\bbr}{\reals}
\newcommand{\bbn}{{\mathbb N}}

\newcommand{\X}{{\mathcal{X}}}
\newcommand{\Y}{{\mathcal{Y}}}

\newcommand{\ellp}{{\ell^{\prime}}}

\newcommand{\one}{{\bf 1}}

\newcommand{\psiinv}{\psi^{\leftarrow}}

\begin{document}

\bibliographystyle{abbrv}

\title[FCLT for subgraph counting processes]
{Functional Central Limit Theorem for Subgraph Counting Processes}
\author{Takashi Owada}
\address{Faculty of Electrical Engineering\\
Technion-Israel Institute of Technology \\
Haifa, 32000, Israel}
\email{takashiowada@ee.technion.ac.il}

\thanks{This research was supported by funding from the European Research Council under the European Union's
Seventh Framework Programme (FP/2007-2013) / ERC Grant Agreement n. 320422.}

\keywords{Extreme value theory, functional central limit theorem, geometric graph, regular variation, von-Mises function. \vspace{.5ex}}

\begin{abstract}
The objective of this study is to investigate the limiting
behavior of a subgraph counting process. The subgraph counting
process we consider counts the number of subgraphs having a
specific shape that exist outside an expanding ball as the sample
size increases. As underlying laws, we consider distributions with either a regularly varying tail or an
exponentially decaying tail. In both cases, the nature of the
resulting functional central limit theorem differs according to the
speed at which the ball expands. More specifically, the normalizations in the central limit theorems and the properties of the limiting Gaussian processes are all determined by whether or not
an expanding ball covers a region - called a weak core - in which
the random points are highly densely scattered and form a giant
geometric graph.
\end{abstract}

\maketitle

\section{Introduction} \label{sec:intro}

The history of random geometric graphs started with Gilbert's
1961 study (\cite{gilbert:1961}) and, since then, it has received
much attention both in theory and applications. More formally,
given a finite set $\mathcal{X} \subset \bbr^d$ and a real number
$r>0$, the geometric graph $G(\mathcal{X},r)$ is defined as an
undirected graph with vertex set $\mathcal{X}$ and edges $[x,y]$
for all pairs $x,y\in\X$ for which $\|x-y\| \leq r$. The theory of
geometric graphs has  been applied mainly in large communication
network analysis, in which the connectivity of network agents
strongly depends on the distance between them; see
\cite{chen:jia:2001}, \cite{stojmenovic:seddigh:zunic:2002}, and
Chapter 3 of \cite{hekmat:2006}. On the purely theoretical side of
 random geometric graphs, the monograph
\cite{penrose:2003} is probably the best known resource. It covers a wide range of topics, such as the asymptotics of the
number of subgraphs with a specific shape, the vertex degree, the
clique number, the formation of a giant component, etc. From among
these interesting subjects, the present study focuses on
constructing the functional central limit theorem (FCLT) for the number of
subgraphs isomorphic to a predefined connected graph $\Gamma$ of
finite vertices.

A typical setup in \cite{penrose:2003} is as follows. Let
$\mathcal X_n$ be a set of random points on $\bbr^d$. Typically, this will be either an $\text{i.i.d.}$ random sample of $n$ points from $f$, or an
inhomogeneous Poisson point process with intensity $nf$, where $f$ is a probability density. We assume that the threshold radius $r_n$
depends on $n$ and decreases to $0$ as $n \to \infty$, but we do
not impose any restrictive assumptions on $f$ except for boundedness. Then, the asymptotic behavior of
the subgraph counts given by
\begin{equation}  \label{e:count.intro1}
G_n := \sum_{\Y \subset \mathcal X_n} \one \bigl\{ G(\Y, r_n) \cong \Gamma \bigr\}\,,
\end{equation}
($\cong$ denotes graph isomorphism, and $\Gamma$ is a fixed
connected graph) splits into three different regimes. First, if
$nr_n^d \to 0$, called the \textit{subcritical} or \textit{sparse}
regime, the distribution of subgraphs isomorphic to $\Gamma$ is
sparse,  and these subgraphs are mostly observed as isolated components. If $nr_n^d \to
\xi \in (0,\infty)$, called the \textit{critical} or
\textit{thermodynamic} regime, for which $r_n$ decreases to $0$ at
a slower rate than the subcritical regime, many of the isolated
subgraphs in $G(\mathcal X_n, r_n)$ become connected to one another.
Finally, if $nr_n^d \to \infty$ (the \textit{supercritical}
regime), the subgraphs are very highly connected and create a
large component.

Historically, the research on the limiting behavior of subgraph
counts of the type \eqref{e:count.intro1} dates back to the
studies of \cite{hafner:1972}, \cite{silverman:brown:1978}, and
\cite{weber:1983}, in all of which mainly the subcritical regime
was treated. Furthermore, \cite{bhattacharya:ghosh:1992} adopted
an approach based on the martingale CLT for $U$-statistics and
proved a CLT under various conditions on $f$ and $r_n$. Relying on
the so-called Stein-Chen method, a set of extensive results for
all three regimes was nicely summarized in Chapter 3 of
\cite{penrose:2003}. Recently, as a higher-dimensional analogue of
a random geometric graph, there has been growing interest in the
asymptotics of the so-called random C\v{e}ch complex. See, for
example, \cite{kahle:2011}, \cite{kahle:meckes:2013}, and
\cite{yogeshwaran:subag:adler:2014}, while
\cite{bobrowski:kahle:2014} provides an elegant review of that
direction.

Somewhat parallel to \eqref{e:count.intro1}, but more
important for the study on the geometric features of extreme sample clouds, is an alternative that we  explore in this
paper. To set this up, we introduce a growing sequence $R_n \to
\infty$ and a threshold radius $t > 0$. The
following quantity, $G_n(t)$ counts the number of subgraphs
in $G(\mathcal X_n, t)$ isomorphic to $\Gamma$ that exist outside
a centered ball in $\bbr^d$ with
radius $R_n$:
\begin{equation}  \label{e:count.intro2}
G_n(t) :=  \sum_{\Y \subset \mathcal X_n } \one \bigl\{ G(\Y, t) \cong \Gamma \bigr\} \times \one \bigl\{ m(\Y) \geq R_n \bigr\}\,,
\end{equation}
where $m(x_1,\dots,x_k) = \min_{1 \leq i \leq k} ||x_i||$, $x_i \in \bbr^d$, and $\| \cdot \|$ is the usual Euclidean norm.

From the viewpoint of extreme value theory (EVT), it is important to investigate limit theorems for $G_n(t)$.  Indeed, over the last decade or so there have been numerous papers treating geometric descriptions of multivariate extremes, among them \cite{balkema:embrechts:2007}, \cite{balkema:embrechts:nolde:2010}, and \cite{balkema:embrechts:nolde:2013}.  In particular, Poisson limits of point processes possessing a U-statistic structure were investigated by \cite{dabrowski:dehling:mikosch:sharipov:2002} and \cite{schulte:thale:2012}, the latter also treating  a number of examples in stochastic geometry. The main references for EVT are \cite{embrechts:kluppelberg:mikosch:1997},
\cite{resnick:1987}, and \cite{dehaan:ferreira:2006}.

The asymptotic behavior of \eqref{e:count.intro2} has been partially explored in \cite{owada:adler:2015}, where a growing sequence $R_n$ is taken in such a way that \eqref{e:count.intro2} has Poisson limits as $n\to\infty$. The main contribution in \cite{owada:adler:2015} is the discovery of a certain layered structure consisting of a collection of ``rings" around the origin with each ring containing extreme random points which exhibit different geometric and topological behavior. The object of the current study is to develop a fuller description of this ring-like structure, at least in a geometric graph model, by establishing a variety of FCLTs which describe geometric graph formation between the rings. 

By  construction, the subgraph counts \eqref{e:count.intro2} can be viewed as generating a stochastic process in the parameter $t\geq0$, while a process-level extension in \eqref{e:count.intro1} is much less obvious. Then, while \eqref{e:count.intro2} captures the dynamic evolution of geometric graphs as $t$ varies,  \eqref{e:count.intro1} only describes the static geometry. Thus, the limits in the FCLT for \eqref{e:count.intro2} are intrinsically Gaussian processes, rather than one-dimensional Gaussian distributions. 

One of the main results of this paper is that the limiting Gaussian processes can be classified into three distinct categories, according to how rapidly $R_n$ grows. 
The most important condition for this classification is whether or not a ball centered at the origin with radius $R_n$, denoted by $B(0, R_n)$, asymptotically covers a \textit{weak core}. Weak cores are balls, centered at the origin with growing radii as $n$ increases, in which the random points are densely scattered and form a highly connected geometric graph. This notion, along with the related notion of a \textit{core}, play a crucial role for the classification of the limiting Gaussian processes. Indeed, if $B(0, R_n)$ grows so that it asymptotically covers a weak core, then the geometric graph outside $B(0, R_n)$ is ``sparse" with many small disconnected components.  In this case, the limit is denoted as the difference between two time-changed Brownian motions. In contrast, if $B(0,R_n)$ is asymptotically covered by a weak core, the geometric graph in the area between the outside of $B(0,R_n)$ and inside of a weak core becomes ``dense", and, accordingly, the limit becomes a degenerate Gaussian process with deterministic sample paths. 
Finally if $B(0,R_n)$ coincides with a weak core, then the limiting Gaussian process possesses more complicated structure and are even non-self-similar. 

We want to emphasize that the nature of the FCLT depends not only on the growth rate of $R_n$ but also the tail property of $f$. This is in complete
contrast to \eqref{e:count.intro1}, because, as seen in Chapter 3
of \cite{penrose:2003}, the proper normalization, limiting
Gaussian distribution, etc. of the CLT are all robust to whether
$f$ has a heavy or a light tail. In this paper, we particularly
deal with the distributions of regularly varying tails and
(sub)exponential tails. However, we are not basically concerned
with any distribution with a superexponential tail, e.g., a
multivariate normal distribution. The details of the FCLT in that case remain for a future study.

The remainder of the paper is organized as follows. First, in
Section 2 we provide a formal definition of the subgraph counting
process. Section 3 gives an overview of what was shown in the previous work \cite{owada:adler:2015} and what will be shown in this paper. Subsequently, in Section 4 we focus on the case in which
the underlying density has a regularly varying tail, including
power-law tails, and prove the required FCLT. We also investigate the
properties of the limiting Gaussian processes, in particular, in
terms of self-similarity and sample path continuity. In
Section 5, we do the same when the underlying density has an
exponentially decaying tail. To distinguish densities via their
tail properties, we need basic tools in EVT. In essence, the properties of the limiting Gaussian processes are determined by how rapidly $R_n$ grows to infinity, as well as how rapidly the tail of $f$ decays. Finally, Section 6 carefully examines both cores and weak cores for a large class of
densities.

Before commencing the main body of the paper, we remark that all the random points in this paper are assumed to be
generated by an inhomogeneous Poisson point process on $\bbr^d$
with intensity $nf$. In our opinion, the FCLT in the main theorem
can be carried over to a usual $\text{i.i.d.}$ random sample setup by a
standard ``de-Poissonization" argument; see Section 2.5 in
\cite{penrose:2003}. This is, however, a little more technical and
challenging, and therefore, we decided to concentrate on the
simpler setup of an inhomogeneous Poisson point process. Furthermore we  consider only spherically symmetric distributions.
Although the spherical symmetry assumption is far from being
crucial, we adopt it to avoid unnecessary technicalities.
%\vspace{10pt}

\section{Subgraph Counting Process}

Let $(X_i, \, i \geq 1)$ be $\text{i.i.d.}$ $\bbr^d$-valued
random variables with spherically symmetric probability density
$f$. Given a Poisson random variable $N_n$ with mean $n$,
independent of $(X_i, \, i \geq 1)$, denote by $\mathcal{P}_n = \{
X_1, X_2, \dots, X_{N_n} \}$ a Poisson point process with
$|\mathcal{P}_n| := N_n$. We choose a positive integer $k$, which
remains fixed hereafter. We take $k \geq 2$, unless otherwise
stated, because many of the functions and objects to follow are degenerate in the case of $k=1$. 

Let $\Gamma$ be a fixed connected graph of $k$
vertices and $G$ represent a geometric graph;  $\cong$ denotes 
graph isomorphism. 
We define
$$
h(x_1, \dots, x_k) := \one \bigl\{ G \bigl( \{ x_1, \dots, x_k \}, \, 1 \bigr) \cong \Gamma \bigr\}\,, \ \ x_1,\dots,x_k \in \bbr^d\,.
$$
Next, we define a collection of indicators $(h_t, \, t
\geq 0)$ by
\begin{equation}  \label{e:geo.graph.dyna}
h_t(x_1,\dots,x_k) := h(x_1/t,\dots,x_k/t) = \one \bigl\{ G \bigl( \{ x_1, \dots, x_k \}, \, t \bigr) \cong \Gamma \bigr\}\,,
\end{equation}
from which one can capture the manner in which a geometric
graph dynamically evolves as the threshold radius $t$ varies. Note, in particular, that $h_1(x_1,\dots,x_k) = h(x_1,\dots,x_k)$.

Clearly $h_t$ is shift invariant: 
\begin{align}
h_t(x_1,\dots,x_k) &= h_t(x_1+y,\dots,x_k+y)\,,\ \ \ x_1,\dots,x_k,y \in \bbr^d\,, \label{e:location.inv}
\end{align}
and, further,
\begin{equation}  \label{e:close.enough}
h_t(0,x_1,\dots,x_{k-1}) = 0 \ \ \text{if } ||x_i|| > kt \ \text{for some } i=1,\dots,k-1\,.
\end{equation}
The latter condition implies that $h_t(x_1,\dots,x_k)=1$ only when all the points
$x_1,\dots,x_k$ are close enough to each other. 

Moreover $h_t$ can be decomposed as follows. Suppose that $\Gamma$
has $k$ vertices and $j$ edges for some $j \in \bigl\{
k-1,\dots,k(k-1)/2 \bigr\}$. Letting $A_\ell$ be a set of
connected graphs of $k$ vertices and $\ell$ edges (up to graph
isomorphism), define for $x_1,\dots,x_k \in \bbr^d$,
\begin{align*}
h_t^+ (x_1,\dots,x_k) &:= h_t(x_1,\dots,x_k) + \sum_{\ell=j+1}^{k(k-1)/2} \sum_{\Gamma^{\prime} \in A_\ell} \one \bigl\{ G\bigl(\{ x_1,\dots,x_k \},t\bigr) \cong \Gamma^{\prime} \big\}\,, \\
h_t^- (x_1,\dots,x_k) &:= \sum_{\ell=j+1}^{k(k-1)/2} \sum_{\Gamma^{\prime} \in A_\ell} \one \bigl\{ G\bigl(\{  x_1,\dots,x_k\},t\bigr) \cong \Gamma^{\prime} \big\}\,.
\end{align*}
Note that $h_t^+(x_1,\dots,x_k)=1$ if and only if a geometric graph $G \bigl( \{ x_1,\dots,x_k \}, t
\bigr)$ either coincides with $\Gamma$ (up to graph isomorphism) or has
more than $j$ edges, while $h_t^-(x_1,\dots,x_k) = 1$ only when $G \bigl( \{ x_1,\dots,x_k \}, t \bigr)$
has more than $j$ edges. It is then elementary to check that
$h_t^{\pm}$ are both indicators, taking values $0$ or $1$, and satisfying, for all $x_1,\dots,x_k \in \bbr^d$ and $0 \leq s \leq t$,
\begin{align}
h_t (&x_1,\dots,x_k) = h_t^+ (x_1,\dots,x_k) - h_t^- (x_1,\dots,x_k)\,, \label{e:ind.decomp1} \\
&h_s^{+} (x_1,\dots,x_k) \leq h_t^{+} (x_1,\dots,x_k)\,, \label{e:ind.increase+}  \\
&h_s^{-} (x_1,\dots,x_k) \leq h_t^{-} (x_1,\dots,x_k)\,.  \notag  \\
&h_t^{\pm}(0,x_1,\dots,x_{k-1}) = 0 \ \ \text{if } ||x_i|| > kt \ \text{for some } i=1,\dots,k-1\,.  \label{e:close.enough.decomp.dyna}
\end{align}
In addition, 
since $h_t$ is an indicator, it is always the case that
$$
h_t^-(x_1,\dots,x_k) \leq h_t^+(x_1,\dots,x_k)\,.
$$

The objective of this study is
to establish a functional central limit theorem (FCLT) of the
\textit{subgraph counting process} defined by
\begin{equation}  \label{e:subgraph.count}
G_n(t) := \sum_{\Y \subset \mathcal{P}_n} h_t(\Y)\, \one \bigl\{ m(\Y) \geq R_n \bigr\}\,,  \ \ t \geq 0\,,
\end{equation}
where $h_t$ is given in \eqref{e:geo.graph.dyna}, $m(x_1,\dots,x_k) = \min_{1 \leq i \leq k} ||x_i||$, $x_i
\in \bbr^d$, and $(R_n, \, n\geq1)$ is a properly chosen
normalizing sequence. Note that \eqref{e:subgraph.count} counts
the number of subgraphs in $G(\Pn, t)$ isomorphic to $\Gamma$ that
lie completely outside of $B(0,R_n)$. More
concrete definitions of $(R_n)$ are given in the subsequent
sections, where the sequence is shown to be dependent on the tail decay rate of $f$.

\section{Annuli Structure}  \label{s:annuli}

The objective of this short section is to clarify what is already known and
what is new in this paper. Without any real loss of generality, we will do this via two simple examples,  one of
which treats a power-law density and the other a density with a (sub)exponential tail. 
Before this, however, we introduce two important notions. \begin{definition}(\cite{adler:bobrowski:weinberger:2014})  \label{def.core}
Given an inhomogeneous Poisson
point process $\Pn$ in $\bbr^d$ with a spherically symmetric density $f$, a
centered ball $B(0,R_n)$, with $R_n \to \infty$, is called a
\textit{core} if
\begin{equation}  \label{e:core.event}
B(0, R_n) \subset \bigcup_{X \in \Pn \cap B(0, R_n)} B(X, 1)\,.
\end{equation}
\end{definition}
In other words, a core is a centered ball in which random
points are densely scattered, so that placing unit balls around
them covers the ball itself. We usually wish to seek the
largest possible value of $R_n$ such that \eqref{e:core.event}
occurs asymptotically with probability $1$. A related notion, the
\textit{weak core}, plays a more decisive role in characterizing
the FCLT proven in this paper. It is shown later that a
weak core is generally larger but close in size to a core of maximum size.
\begin{definition}  \label{def.weak.core}
Let $f$ be a spherically symmetric density on $\bbr^d$ and $e_1 =
(1,0,\dots,0) \in \bbr^d$. A \textit{weak core} is a centered ball
$B(0,R_n^{(w)})$ such that $nf(R_n^{(w)}e_1) \to 1$ as $n \to
\infty$.
\end{definition}

\begin{example}  \label{ex:power.law.tail}
{\rm Consider the power-law density
\begin{equation}  \label{e:simple.pdf.RV}
f(x) = C/\bigl( 1 + ||x||^{\alpha} \bigr)\,, \ \ x \in \bbr^d,
\end{equation}
for some $\alpha > d$ and normalizing constant $C$. Using this density, we  see how random
geometric graphs are formed in all of $\bbr^d$. First, according
to \cite{adler:bobrowski:weinberger:2014}, there exists a sequence
$R_n^{(c)} \sim \text{constant} \times (n/\log n)^{1/\alpha}$, $n
\to \infty$ such that, if $R_n \leq R_n^{(c)}$,
\eqref{e:core.event} occurs asymptotically with probability $1$.
%The authors concluded that this $R_n^{(c)}$ is the radius of amaximum core. 
In addition, as for the radius of a weak core, it
suffices to take $R_n^{(w)} = (Cn)^{1/\alpha}$. Although
$R_n^{(w)}$ grows faster than $R_n^{(c)}$, they are seen to be
``close" to each other in the sense that they have the same
regular variation exponent, $1/\alpha$.

Beyond a weak core, however, the formation of random geometric
graphs drastically varies. In fact, the exterior of a weak core
can be divided into annuli of different radii, at which many
isolated subgraphs of finite vertices are asymptotically placed in
a specific fashion. To be more precise, let us fix connected
graphs $\Gamma_k$ with $k$ vertices for $k=2,3,\dots$ and let
$$
R_{k,n}^{(p)} := \bigl( C n\bigr)^{1/(\alpha - d/k)},
$$
which in turn implies that $R_n^{(w)} \ll \cdots \ll R_{k,n}^{(p)} \ll R_{k-1,n}^{(p)} \ll \cdots \ll R_{2,n}^{(p)}$, and
$$
n^k \bigl( R_{k,n}^{(p)} \bigr)^d f\bigl( R_{k,n}^{(p)} e_1 \bigr)^k \to 1\,, \ \ n \to \infty\,.
$$
Under this circumstance, \cite{owada:adler:2015} considered the subgraph counts given by
\begin{equation}  \label{e:count.intro3}
\sum_{\Y \subset \Pn} \one \bigl\{ G(\Y,t) \cong \Gamma_k \bigr\} \times \one \bigl\{ m(\Y) \geq R_{k,n}^{(p)} \bigr\}\,,
\end{equation}
and showed that \eqref{e:count.intro3} weakly converges to a
Poisson distribution for each fixed $t$. To be more specific on
the geometric side, let Ann$(K,L)$ be an annulus with inner radius
$K$ and outer radius $L$. Then, we have, in an asymptotic sense,
\begin{itemize}
\item Outside $B\bigl(0,R_{2,n}^{(p)}\bigr)$, there are finitely
many graphs isomorphic to $\Gamma_2$, but  none isomorphic to  $\Gamma_3, \Gamma_4, \dots$. \item
Outside $B\bigl(0,R_{3,n}^{(p)}\bigr)$, equivalently inside
Ann$\bigl( R_{3,n}^{(p)}, R_{2,n}^{(p)} \bigr)$, there are
infinitely many graphs isomorphic to $\Gamma_2$ and finitely many graphs isomorphic to $\Gamma_3$,  but none isomorphic to
 $\Gamma_4, \Gamma_5, \dots$.
\end{itemize}
In general,
\begin{itemize}
\item Outside $B\bigl(0,R_{k,n}^{(p)}\bigr)$, equivalently inside
Ann$\bigl( R_{k,n}^{(p)}, R_{k-1,n}^{(p)} \bigr)$, there are
infinitely many graphs isomorphic to $\Gamma_2, \dots, \Gamma_{k-1}$ and finitely many
graphs isomorphic to  $\Gamma_k$,  but none isomorphic to  $\Gamma_{k+1}, \Gamma_{k+2}, \dots$ etc.
\end{itemize}

Section \ref{s:heavy} of the current paper considers the subgraph
counts of the form 
\begin{equation}  \label{e:count.intro4}
\sum_{\Y \subset \Pn} \one \bigl\{ G(\Y,t) \cong \Gamma_k \bigr\} \times \one \bigl\{ m(\Y) \geq R_n \bigr\}\,,
\end{equation}
where $(R_n)$ satisfies
\begin{equation}  \label{e:rate.clt}
n^k R_n^d f(R_ne_1)^k \to \infty\,, \ \ n \to \infty\,,
\end{equation}
in which case, $R_n \ll R_{k,n}^{(p)}$. As a consequence of
\eqref{e:rate.clt}, we may naturally anticipate that a
FCLT governs the asymptotic behavior of \eqref{e:count.intro4}. Since $(R_n)$ satisfying
\eqref{e:rate.clt} shows a slower divergence rate than
$(R_{k,n}^{(p)})$, i.e., $R_n/R_{k,n}^{(p)} \to 0$, we may expect
that infinitely many subgraphs isomorphic to $\Gamma_k$ appear
asymptotically outside $B(0,R_n)$. This in turn implies that,
instead of a Poisson limit theorem, the FCLT governs the limiting
behavior of the subgraph counting process.

%In view of the layered
%structure described above, a Poisson limit theorem no
%longer governs the behavior of \eqref{e:count.intro4};
%Alternatively, one needs to establish a central limit
%theorem. 

As the analog of the setup for \eqref{e:count.intro1}, when deriving an FCLT, the behavior
of \eqref{e:count.intro4} splits into three different regimes:
$$
(i)\ nf(R_ne_1) \to 0\,, \ \ (ii)\ nf(R_ne_1) \to \xi \in (0,\infty)\,, \ \ (iii)\ nf(R_ne_1) \to \infty\,.
$$
Specifically, if $nf(R_ne_1) \to 0$ (i.e., $B(0,R_n)$ contains a
weak core), many isolated components of subgraphs isomorphic to
$\Gamma_k$ are distributed outside $B(0,R_n)$. If
$nf(R_ne_1) \to \xi \in (0,\infty)$ (i.e., $B(0,R_n)$ agrees with a weak core), the subgraphs isomorphic
to $\Gamma_k$ outside $B(0,R_n)$ begin to be connected to one another. In particular, observing that $\lim_{k\to\infty} R_{k,n}^{(p)} = R_n^{(w)}$ for all $n$, we see that 
\begin{itemize}
\item Outside of $B(0, R_n^{(w)})$, there are infinitely many graphs isomorphic to $\Gamma_j$ for every $j=2,3,\dots$.
\end{itemize}
If $nf(R_ne_1) \to \infty$ (i.e.,  $B(0,R_n)$ is contained in a weak core), the subgraphs isomorphic
to $\Gamma_k$ outside $B(0,R_n)$ are further increasingly connected and form a large component. 

In Section \ref{s:heavy}, we will see that the nature of the FCLT, including the normalizing constants and the properties of the limiting Gaussian processes, differs according to which regime one considers. Combing the results on the FCLT and the Poissonian results in \cite{owada:adler:2015}, we obtain a complete picture of the annuli structure formed by heavy tailed random variables. 

\begin{figure}[!t]
  \label{f:layered.structure}
\includegraphics[width=13.5cm]{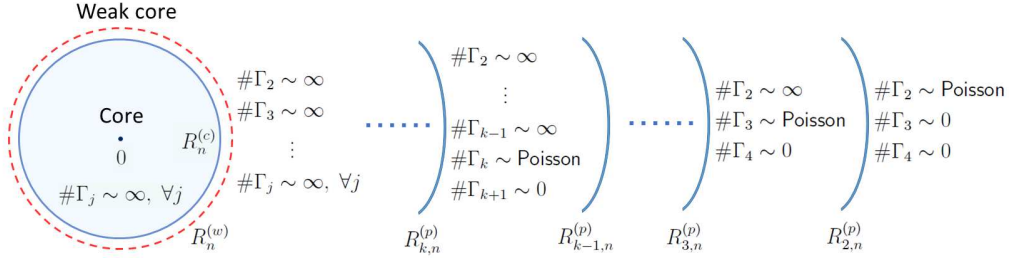}
\vspace{-10pt}
\caption{{\footnotesize Layered structure of random geometric graphs. For the
density \eqref{e:simple.pdf.RV}, $R_n^{(c)}$ and $R_n^{(w)}$ are
regularly varying sequences with exponent $\alpha^{-1}$.
$R_{k,n}^{(p)}$ is also a regularly varying sequence with exponent $(\alpha -
d/k)^{-1}$.  We study the FCLT for \eqref{e:count.intro4}
in three different regimes, i.e., $(i)\ nf(R_ne_1) \to 0$, $(ii)\
nf(R_ne_1) \to \xi \in (0,\infty)$, and $(iii)\ nf(R_ne_1) \to
\infty$. In relation to other radii, they are respectively
equivalent to $(i)\ R_n^{(w)} \ll R_n \ll R_{k,n}^{(p)}$, $(ii)\
R_n \sim R_n^{(w)}$, and $(iii)\ R_n \ll R_n^{(w)}$.}}
\end{figure}
}
\end{example}
\begin{example}  \label{ex:subexp.tail}
{\rm Next, we turn to a density with a (sub)exponential tail
$$
f(x) = C e^{-||x||^{\tau}/\tau}, \ \ x \in \bbr^d\,, \ 0 < \tau \leq 1\,.
$$
for which the radius of a maximum core is given by
$$
R_n^{(c)} = \bigl(\tau \log n - \tau \log \log (\tau \log n)^{1/\tau} + \text{constant} \bigr)^{1/\tau};
$$
see \cite{adler:bobrowski:weinberger:2014} and \cite{owada:adler:2015}. Obviously, one can
take $R_n^{(w)} = \bigl(\tau \log n + \tau \log C \bigr)^{1/\tau}$. As in the previous example, the exterior of a weak core is characterized by the
same kind of layer structure, for which the description in Figure 1 applies, except for the change in the values of $R_{k,n}^{(p)}$. Letting
$$
R_{k,n}^{(p)} = \bigl( \tau \log n + k^{-1} (d-\tau) \log (\tau \log n) + \tau \log C \bigr)^{1/\tau},
$$
we have, in an asymptotic sense,
$R_n^{(w)} \ll \cdots \ll R_{k,n}^{(p)} \ll R_{k-1,n}^{(p)} \ll \cdots \ll R_{2,n}^{(p)}$, and
$$
n^k \bigl( R_{k,n}^{(p)} \bigr)^{d-\tau} f\bigl( R_{k,n}^{(p)} e_1 \bigr)^k \to 1\,, \ \ n \to \infty\,.
$$
Then, it was shown in \cite{owada:adler:2015} that
\eqref{e:count.intro3} converges weakly to a Poisson distribution
for each fixed $t$.  

In Section \ref{s:light} of this paper, taking $(R_n)$ such that
$n^k R_{n}^{d-\tau} f( R_{n} e_1 )^k \to \infty$, we establish a
FCLT for the subgraph counting process \eqref{e:subgraph.count}. To this end, our
argument has to be split, once again, into the three different
regimes:
$$
(i)\ nf(R_ne_1) \to 0\,, \ \ (ii)\ nf(R_ne_1) \to \xi \in (0,\infty)\,, \ \ (iii)\ nf(R_ne_1) \to \infty\,.
$$
}
As in the last example, three different Gaussian limits may appear depending on  the regime. This completes the full description of the annuli structure formed by random variables with an exponentially decaying tail, when combined with the Poisson limit theorems in \cite{owada:adler:2015}.
\end{example}
%\vspace{10pt}

\section{Heavy Tail Case}  \label{s:heavy}

\subsection{The Setup}
In this section, we explore the case in which the underlying
density $f$ on $\bbr^d$ has a heavy tail under a more general setup than that in Example \ref{ex:power.law.tail}. Let $S_{d-1}$ be a $(d-1)$-dimensional unit sphere in $\bbr^d$. We
assume that the density has a regularly varying tail (at infinity)
in the sense that for any $\theta \in S_{d-1}$ (equivalently, for
some $\theta \in S_{d-1}$ because of the spherical symmetry of
$f$), and for some $\alpha > d$,
$$
\lim_{r \to \infty} \frac{f(rt \theta)}{f(r \theta)} = t^{-\alpha} \ \ \text{for every } t>0\,.
$$
Denoting by $RV_{-\alpha}$ a collection of regularly varying
functions (at infinity) of exponent $-\alpha$, the above is
written as
\begin{equation}  \label{e:RV.tail}
f \in RV_{-\alpha}\,.
\end{equation}
Clearly, a power-law density in Example \ref{ex:power.law.tail} satisfies \eqref{e:RV.tail}.
Let $k \geq 2$ be an integer that remains fixed throughout this
section. We remark that many of the functions and objects are
dependent on $k$, but the dependence may not be stipulated by
subscripts (or superscripts). Choosing the sequence $R_n
\to \infty$ so that
\begin{equation}  \label{e:normalizing.heavy}
n^k R_n^d f(R_n e_1)^k \to \infty \ \ \text{as } n \to \infty\,,
\end{equation}
we consider the subgraph counting process given in
\eqref{e:subgraph.count}, whose behavior is, as argued in Example \ref{ex:power.law.tail}, expected to be governed by a FCLT.

The scaling constants for the FCLT, denoted by $\tau_n$, are shown
to depend on the limit value of $nf(R_ne_1)$ as $n \to \infty$.
More precisely, we take
\begin{equation}  \label{e:tau.heavy}
\tau_n := \begin{cases} n^k R_n^d f(R_ne_1)^k & \text{if } nf(R_ne_1) \to 0\,, \\ R_n^d & \text{if } nf(R_ne_1) \to \xi \in (0,\infty)\,, \\ n^{2k-1} R_n^d f(R_ne_1)^{2k-1} & \text{if } nf(R_ne_1) \to \infty\,.
\end{cases}
\end{equation}
The reason for which we need three different normalizations is
deeply related to the connectivity of a random geometric graph. To
explain this, we need the notion of a \textit{weak core}; see Definition \ref{def.weak.core} for the formal definition. The main point is that the density
of random points between the outside and inside of a weak core
is completely different. In essence, random points inside a weak
core are highly densely scattered, and the corresponding random
geometric graph forms a single giant component. Beyond a weak
core, however, random points are distributed less densely, and as
a result, we observe many isolated geometric graphs of smaller
size. This disparity between the outside and inside of a weak
core requires different normalizations in $(\tau_n)$. In Section
\ref{s:connect.core},  a more detailed study in this direction is
presented. \vspace{5pt}

\subsection{Limiting Gaussian Processes and the FCLT}  \label{s:limit.heavy}

We introduce a family of Gaussian processes which function as the building blocks for the limiting Gaussian processes in the FCLT. For $\ell = 1,\dots, k$, let
$$
B_\ell = \frac{s_{d-1}}{\ell! \bigl( (k-\ell)! \bigr)^2 \bigl( \alpha(2k-\ell) - d \bigr)}\,,
$$
where $s_{d-1}$ is a surface area of the $(d-1)$-dimensional unit sphere in $\bbr^d$. 

For $\ell = 2,\dots,k$, write $\lambda_\ell$ for the Lebesgue measure on $(\bbr^d)^{\ell-1}$, and denote by $G_\ell$ a \textit{Gaussian $B_\ell \lambda_\ell$-noise}, such that 
$$
G_\ell(A) \sim \mathcal N \bigl( 0, B_\ell \lambda_\ell (A) \bigr)
$$
for measurable sets $A \subset (\bbr^d)^{\ell-1}$ with $\lambda_\ell (A) < \infty$, and if $A \cap B= \emptyset$, then $G_\ell(A)$ and $G_\ell (B)$ are independent. For $\ell = 1$, we define $G_1$  as a Gaussian random variable with zero mean and variance $B_1$. We assume that $G_1, \dots, G_k$ are independent. 

For $\ell = 2,\dots, k-1$, we define Gaussian processes $\BV_\ell  = \bigl(V_\ell (t), \, t \geq 0\bigr)$ by 
$$
V_\ell (t) := \int_{(\bbr^d)^{\ell-1}} \int_{(\bbr^d)^{k-\ell}} h_t (0, \by, \bz) \, d\bz\,  G_\ell (d\by), \ \ t\geq0. 
$$
In addition, if $\ell = k$, define 
$$
V_k(t) := \int_{(\bbr^d)^{k-1}} h_t(0,\by) G_k (d\by),
$$
and if $\ell = 1$, set 
$$
V_1(t) := \int_{(\bbr^d)^{k-1}} h_t(0,\bz)\, d\bz\, G_1 = t^{d(k-1)} \int_{(\bbr^d)^{k-1}} h(0,\bz)\, d\bz\, G_1.
$$
Note that $\BV_1$ is a degenerate Gaussian process with deterministic sample paths. These processes later turn out to be the
building blocks of the weak limits in the main theorem.

The covariance function of the process $\BV_\ell$ is given by 
\begin{align}  
L_\ell(t,s) &:= \E \bigl\{ V_\ell(t) V_\ell(s) \bigr\} \label{e:cov.comp} \\
&= B_\ell \int_{(\bbr^d)^{\ell-1}} \hspace{-10pt} d\by \int_{(\bbr^d)^{k-\ell}}\hspace{-10pt} d\bz_2 \int_{(\bbr^d)^{k-\ell}}  \hspace{-10pt} d\bz_1\, h_t(0,\by,\bz_1)\, h_s(0,\by,\bz_2)\,, \ \ \ t,s \geq 0 \notag
\end{align}
(if $\ell =k$, we take $\bz_i = \emptyset$, $i=1,2$, and if $\ell = 1$, we set $\by = \emptyset$).

Using the decomposition \eqref{e:ind.decomp1}, we can express $\BV_\ell$ as the difference between two Gaussian processes; that is, for $\ell = 2,\dots,k-1$,
\begin{align*}
V_\ell (t) &= \int_{(\bbr^d)^{\ell-1}} \int_{(\bbr^d)^{k-\ell}} \hspace{-5pt}h_t^+ (0, \by, \bz) \, d\bz\,  G_\ell (d\by) - \int_{(\bbr^d)^{\ell-1}} \int_{(\bbr^d)^{k-\ell}} \hspace{-5pt}h_t^- (0, \by, \bz) \, d\bz\,  G_\ell (d\by) \\
&:= V_\ell^+(t) - V_\ell^-(t).
\end{align*}
The same decomposition is feasible in an analogous manner for $\BV_1$ and $\BV_k$. 

The following proposition shows that the processes $\BV_k^{+}$ and $\BV_k^-$ can be represented as a time-changed Brownian motion. 
\begin{proposition}  \label{p:limit.heavy1}
The process $\BV_k^+$ can be expressed as 
$$
\bigl( V_k^+(t), \, t\geq 0 \bigr) \stackrel{d}{=} \Bigl( B \bigl( K_k^+ \, t^{d(k-1)} \bigr), \, t \geq 0 \Bigr),
$$
where $B$ is the standard Brownian motion, and $K_k^+ := B_k \int_{(\bbr^d)^{k-1}} h^+(0,\by) d\by$. \\
Replacing $K_k^+$ with $K_k^- := B_k \int_{(\bbr^d)^{k-1}} h^-(0,\by) d\by$, we obtain the same statement for $\BV_k^-$.
\end{proposition}
\begin{proof}
It is enough to verify that the covariance functions on both sides coincide. It follows from \eqref{e:ind.increase+} that for $0 \leq s \leq t$, 
\begin{align*}
\E \bigl\{ V_k^+(t) V_k^+(s)  \bigr\} &= B_k \int_{(\bbr^d)^{k-1}} h_t^+(0,\by)\, h_s^+(0,\by) d\by \\
&= s^{d(k-1)} K_k^+ \\
&= \E \bigl\{ B(K_k^+\, t^{d(k-1)})B(K_k^+\, s^{d(k-1)}) \bigr\}.
\end{align*}
\end{proof}

We also claim that the
process $\BV_\ell$ is self-similar and has a.s. $\hspace{-8pt}$
H\"{o}lder continuous sample paths. Recall that a stochastic
process $\bigl(X(t), \, t \geq 0  \bigr)$ is said to be
self-similar with exponent $H$ if
$$
\bigl(X(ct_i), \, i=1,\dots,k  \bigr) \stackrel{d}{=} \bigl(c^HX(t_i), \, i=1,\dots,k  \bigr)
$$
for any $c>0$, $t_1,\dots,t_k \geq 0$, and $k \geq 1$.

\begin{proposition}  \label{p:self.similar.heavy}
$(i)$ For $\ell=1,\dots,k$, the process $\BV_\ell$ is self similar with exponent $H=d(2k-\ell-1)/2$. \\
\noindent $(ii)$ For $\ell  =1,\dots,k$ and every $T>0$, $\bigl( V_\ell(t), \, 0\leq t
\leq T \bigr)$ has a modification, the sample paths of which are
H\"{o}lder continuous of any order in $[0,1/2)$.
\end{proposition}
\begin{proof}
We can immediately prove $(i)$ by the scaling property 
$$
L_\ell(ct,cs) = c^{d(2k-\ell-1)}L_\ell(t,s)\,, \ \ \ t,s\geq0\,, \ c>0\,.
$$

As for $(ii)$, the statement is obvious for $\ell = 1$ or $\ell =k$; therefore, we take $\ell \in \{ 2,\dots,k-1 \}$. By Gaussianity,
\begin{equation}  \label{e:2m.moment}
\E \Bigl\{ \bigl( V_\ell(t) - V_\ell(s) \bigr)^{2m}\Bigr\} = \prod_{i=1}^m(2i-1)\, \Bigl( \E \Bigl\{ \bigl( V_\ell(t) - V_\ell(s) \bigr)^2 \Bigr\}\Bigr)^m, \ \ m=1,2,\dots
\end{equation}
We now show that there exists a constant $C>0$, which depends on $T$,
such that
\begin{equation}  \label{e:est.2nd.moment}
\E \Bigl\{ \bigl(V_\ell(t) - V_\ell(s) \bigr)^2\Bigr\} \leq C(t-s) \ \ \text{for all } 0 \leq s \leq t \leq T\,.
\end{equation}
By virtue of the decomposition $\BV_\ell = \BV_\ell^+ -
\BV_\ell^-$, showing \eqref{e:est.2nd.moment} for each of
$\BV_\ell^+$ and $\BV_\ell^-$ suffices. We handle $\BV_\ell^+$
only, since $\BV_\ell^-$ can be treated in the same manner. 
We have
\begin{align*}
\E \Bigl\{ \bigl( V_\ell^+(t) - V_\ell^+(s) \bigr)^2\Bigr\}  &= B_\ell \int_{(\bbr^d)^{\ell-1}} \hspace{-10pt} d\by \int_{(\bbr^d)^{k-\ell}}  \hspace{-10pt} d\bz_2 \int_{(\bbr^d)^{k-\ell}}  \hspace{-10pt} d\bz_1 \bigl\{ h_t^+(0,\by,\bz_1) - h_s^+(0,\by,\bz_1) \bigr\} \\
&\quad \times \bigl\{ h_t^+(0,\by,\bz_2) - h_s^+(0,\by,\bz_2) \bigr\}\,.
\end{align*}
Because of \eqref{e:close.enough.decomp.dyna}, the above integral is not
altered if the integral domain is restricted to $(\bbr^d)^{\ell-1}
\times (\bbr^d)^{k-\ell} \times \bigl( B(0,kT) \bigr)^{k-\ell}$.
In addition, by \eqref{e:ind.increase+}, there exist constants
$C_1, C_2>0$, both depending on $T$, such that
\begin{align*}
\E \Bigl\{ \bigl( V_\ell^+(t) - V_\ell^+(s) \bigr)^2\Bigr\} &\leq C_1 \int_{(\bbr^d)^{\ell-1}} \int_{(\bbr^d)^{k-\ell}}  \bigl\{ h_t^+(0,\by,\bz) - h_s^+(0,\by,\bz) \bigr\} d\bz d\by \\
&= C_1 \int_{(\bbr^d)^{\ell-1}} \int_{(\bbr^d)^{k-\ell}}  h^+ (0,\by,\bz)  d\bz d\by\, \bigl( t^{d(k-1)} - s^{d(k-1)} \bigr) \\
&\leq C_2 \int_{(\bbr^d)^{\ell-1}} \int_{(\bbr^d)^{k-\ell}}  h^+ (0,\by,\bz) d\bz d\by\, (t-s) \ \ \text{for all } 0 \leq s \leq t \leq T\,,
\end{align*}
which verifies \eqref{e:est.2nd.moment}. \\
Combining \eqref{e:2m.moment} and \eqref{e:est.2nd.moment}, we have that for some $C_3>0$,
$$
\E \Bigl\{ \bigl(V_\ell(t) - V_\ell(s)\bigr)^{2m} \Bigr\} \leq C_3(t-s)^m \ \ \text{for all } 0 \leq s \leq t \leq T\,.
$$
It now follows from the Kolmogorov continuity theorem that there
exists a modification of $\bigl( V_\ell(t), \, 0 \leq t \leq T
\bigr)$, the sample paths of which are H\"{o}lder continuous of
any order in $\bigl[0,(m-1)/(2m) \bigr)$. Since $m$ is arbitrary,
we are done by letting $m \to \infty$.
\end{proof}

We are now ready to state the FCLT for the subgraph counting process, suitably scaled and centered in such a way that
$$
X_n(t) = \tau_n^{-1/2} \bigl( G_n(t) - \E \{ G_n(t) \} \bigr)\,, \ \ t \geq 0\,.
$$
In the following, $\Rightarrow$ denotes weak convergence. All weak
convergence hereafter are in the space
$\mathcal{D}[0,\infty)$ of right-continuous functions with left
limits. The proof of the theorem is deferred to Section
\ref{s:proof.heavy}.
\begin{theorem}  \label{t:main.heavy}
$(i)$ If $nf(R_ne_1) \to 0$ as $n \to \infty$, then
$$
\bigl( X_n(t), \, t\geq 0 \bigr) \Rightarrow \bigl( V_k(t), \, t\geq 0 \bigr) \ \ \text{in } \mathcal{D}[0,\infty)\,.
$$
\noindent $(ii)$ If $nf(R_ne_1) \to \xi \in (0,\infty)$ as $n \to \infty$, then
$$
\bigl( X_n(t), \, t\geq 0 \bigr) \Rightarrow \left( \sum_{\ell=1}^k \xi^{2k-\ell} V_\ell(t), \, t\geq 0 \right)  \ \ \text{in } \mathcal{D}[0,\infty)\,.
$$
\noindent $(iii)$ If $nf(R_ne_1) \to \infty$ as $n \to \infty$, then
$$
\bigl( X_n(t), \, t\geq 0 \bigr) \Rightarrow \bigl( V_1(t), \, t\geq 0 \bigr) \ \ \text{in } \mathcal{D}[0,\infty)\,.
$$
\end{theorem}
The processes $\BV_1, \dots,
\BV_k$ can be viewed as the building blocks of the limiting
Gaussian processes; however, how many and which ones contribute to the
limit depends on whether the ball $B(0,R_n)$ covers a weak core or
not. If $B(0,R_n)$ covers a weak core, equivalently, $nf(R_ne_1)
\to 0$, then $\BV_k$ is the only process remaining in the limit. Although, as seen in Proposition \ref{p:limit.heavy1}, $\BV_k$ is generally represented as the difference in two time-changed Brownian motions, it can be denoted as a
\textit{single} time-changed Brownian motion when $h_t$ is
increasing in $t$, i.e., $h_s(\Y) \leq h_t(\Y)$ for all $0 \leq s
\leq t$, $\Y \in (\bbr^d)^k$. This is the case when $\Gamma$ is a complete graph, in which case the negative part
$h_t^{-}$ is identically zero.
In contrast, the process $\BV_1$, a degenerate Gaussian process with deterministic sample paths, only appears in the limit when $B(0,R_n)$ is
contained in a weak core, i.e., $nf(R_ne_1) \to \infty$. Finally,
if $B(0,R_n)$ agrees with a weak core (up to multiplicative
constants), all of the processes $\BV_1, \dots, \BV_k$ contribute
to the limit. Interestingly, only in this case, do the weak
limits become non-self-similar. 

\section{Exponentially Decaying Tail Case}  \label{s:light}

\subsection{The Setup}
This section develops the FCLT of the subgraph counting process
suitably scaled and centered, when the underlying density on
$\bbr^d$ possesses an exponentially decaying tail. Typically, in
the spirit of extreme value theory, a class of multivariate
densities with exponentially decaying tails can be formulated by
the so-called \textit{von Mises functions}. See for example,
\cite{balkema:embrechts:2004} and \cite{balkema:embrechts:2007}.
In particular, in the one-dimensional case ($d=1$), the von Mises
function plays a decisive role in the characterization of the
max-domain of attraction of the Gumbel law. See Proposition 1.4 in
\cite{resnick:1987}. We assume that the density $f$ on $\bbr^d$ is
given by
\begin{equation}  \label{e:density.light}
f(x) = L \bigl( ||x|| \bigr) \exp \bigl\{ -\psi \bigl( ||x|| \bigr) \bigr\}\,, \ \ x \in \bbr^d.
\end{equation}
Here, $\psi: \bbr_+ \to \bbr$ is a function of $C^2$-class and is referred to as a von Mises function, so that
\begin{equation}  \label{e:von.Mises}
\psi^{\prime}(z) > 0, \ \ \psi(z) \to \infty, \ \ \bigl( 1/\psi^{\prime} \bigr)^{\prime}(z) \to 0
\end{equation}
as $z \to z_{\infty} \in (0,\infty]$. In this paper, we restrict
ourselves to an unbounded support of the density, i.e.,
$z_{\infty} \equiv \infty$. For notational ease, we introduce the
function $a(z) = 1/\psi^{\prime}(z)$, $z > 0$. Since
$a^{\prime}(z) \to 0$ as $z \to \infty$, the Ces\`{a}ro mean of
$a^{\prime}$ converges as well:
\begin{equation}  \label{e:auxi.slow}
\frac{a(z)}{z} = \frac{1}{z} \int_0^z a^{\prime}(r)dr \to 0\,, \ \ \text{as } z \to \infty\,.
\end{equation}

Suppose that a measurable function $L:\bbr_+ \to \bbr_+$ is \textit{flat} for $a$, that is,
\begin{equation}  \label{e:flat}
\frac{L \bigl( t + a(t) v \bigr)}{L(t)} \to 1 \ \ \text{as } t \to \infty \ \text{uniformly on bounded } v\text{-sets}.
\end{equation}
This condition implies that $L$ behaves as a constant locally in
the tail of $f$, and thus, only $\psi$ plays a dominant role in
the characterization of the tail of $f$. Here, we need to put an
extra technical condition on $L$. Namely, there exist $\gamma \geq
0$, $z_0 > 0$, and $C \geq 1$ such that
\begin{equation}  \label{e:poly.upper}
\frac{L(zt)}{L(z)} \leq C\, t^{\gamma} \ \ \text{for all } t > 1, \, z \geq z_0\,.
\end{equation}

Since $L$ is negligible in the tail of $f$, it seems reasonable to
classify the density \eqref{e:density.light} in terms of the limit
of $a$. If $a(z) \to \infty$ as $z \to \infty$, we say that $f$
belongs to a class of densities with \textit{subexponential} tail, because
the tail of $f$ decays more slowly than that of an exponential
distribution. Conversely, if $a(z) \to 0$ as $z \to \infty$, $f$
is said to have a \textit{superexponential} tail, and if $a(z) \to
c \in (0,\infty)$, we say that $f$ has an exponential tail. To be
more specific about the difference in tail behaviors, let us
consider a slightly more general example than that in Example \ref{ex:subexp.tail}, for which $f(x) = L \bigl( ||x|| \bigr)
\exp \bigl\{ -||x||^{\tau} /\tau \bigr\}$, $\tau > 0$, $x \in \bbr^d$.
Clearly, the parameter $\tau$ is associated with the speed at
which $f$ vanishes in the tail. Observe that $a(z)= z^{1-\tau} \to \infty$ as $z \to \infty$ if $0 < \tau <1$, and
therefore in this case, $f$ has a subexponential tail. If $\tau
> 1$, $a(z)$ decreases to $0$, in which case $f$ has a
superexponential tail.

An important assumption throughout most of this study is that
there exists $c \in (0,\infty]$ such that
\begin{equation}  \label{e:subexp.tail}
a(z) \to c \ \ \text{as } z \to \infty\,.
\end{equation}
In view of the classification described above,
\eqref{e:subexp.tail} eliminates the possibility of densities with superexponential tail. As discovered in
\cite{owada:adler:2015} and
\cite{adler:bobrowski:weinberger:2014}, random points drawn from a
superexponential law hardly form isolated geometric graphs outside
a core, whereas random points coming from a subexponential law do
constitute a layer of isolated geometric graphs outside a core.
Accordingly, it is highly likely that the nature of the FCLT
differs according to whether the underlying density has a
superexponential or a subexponential tail.  The present work
focuses on the (sub)exponential tail case, and more detailed
studies on a superexponential tail case remain for future work.

To realize a more formal set up, let $k \geq 2$ be an integer,
which remains fixed for the remainder of this section; however,
once again, note that many of the functions and objects are
implicitly dependent on $k$. Define the sequence $R_n \to \infty$,
so that
\begin{equation}  \label{e:normalizing.light}
n^k a(R_n) R_n^{d-1} f(R_ne_1)^k \to \infty\,, \ \ \ n \to \infty\,.
\end{equation}
Defining an alternative sequence $R_{k,n}^{(p)} \to \infty$ for which
$$
n^k a\bigl(R_{k,n}^{(p)}\bigr) \bigl(R_{k,n}^{(p)}\bigr)^{d-1} f\bigl(R_{k,n}^{(p)}e_1\bigr)^k \to 1\,, \ \ \ n \to \infty\,,
$$
the subgraph counting process using $R_{k,n}^{(p)}$ is known to
weakly converge to a Poisson distribution; see
\cite{owada:adler:2015}. Since $R_n$ in
\eqref{e:normalizing.light} grows more slowly than
$R_{k,n}^{(p)}$, i.e., $R_n / R_{k,n}^{(p)} \to 0$, we may expect
that an FCLT plays a decisive role in the asymptotic behavior of a
subgraph counting process.

As in the last section, we now want to recall the notion of a
\textit{weak core}. Let $R_n^{(w)} \to \infty$ be a sequence such
that $nf(R_n^{(w)}e_1) \to 1$ as $n \to \infty$. Then, we say that
a ball $B \bigl( 0,R_n^{(w)} \bigr)$ is a weak core. We have to
change, once again, the scaling constants $\tau_n$ of the FCLT,
depending on whether $B(0, R_n)$ covers a weak core or not. More
specifically, we define
\begin{equation}  \label{e:tau.light}
\tau_n := \begin{cases} n^k a(R_n) R_n^{d-1} f(R_ne_1)^k & \text{if } nf(R_ne_1) \to 0\,, \\ a(R_n)R_n^{d-1} & \text{if } nf(R_ne_1) \to \xi \in (0,\infty)\,, \\ n^{2k-1} a(R_n) R_n^{d-1} f(R_ne_1)^{2k-1} & \text{if } nf(R_ne_1) \to \infty\,.
\end{cases}
\end{equation}

\subsection{Limiting Gaussian Processes and the FCLT}  \label{s:limit.proc.light}
The objective of this subsection is to formulate the limiting Gaussian processes and the FCLT. Let
\begin{equation}  \label{e:def.D.ell}
D_\ell = \frac{s_{d-1}}{\ell ! \, \bigl( (k-\ell)! \bigr)^2}\,,  \ \ \ \ell=1,\dots,k,
\end{equation}
and let $H_\ell$ be a Gaussian $\mu_\ell$-noise, where the $\mu_\ell$ for $\ell = 2,\dots,k$, satisfy
\begin{align*}
\mu_\ell(d\rho\, d\by) &= D_\ell\,  e^{-\ell \rho - c^{-1} \sum_{i=1}^{\ell-1} \langle e_1, y_i \rangle}\, \\
&\qquad \times \one \bigl\{ \, \rho + c^{-1} \langle e_1,y_i \rangle \geq 0, \ i = 1,\dots, \ell-1 \,  \bigr\}\, d\rho\, d\by, \ \ \rho \geq 0, \, \by \in (\bbr^d)^{\ell -1}, 
\end{align*}
and 
$$
\mu_1(d\rho) = D_1\, e^{-\rho} d\rho, \ \ \rho \geq 0.
$$
Assume that $H_1, \dots, H_k$ are independent. 

We now define a collection of Gaussian processes  needed for the construction of the limits in the FCLT. For $\ell = 2,\dots, k$,  we define 
\begin{align*}
W_\ell (t) &:= \int_{[0,\infty)\times (\bbr^d)^{\ell -1}} \int_{(\bbr^d)^{k-\ell}} e^{-\sum_{i=1}^{k-\ell} \bigl( \rho + c^{-1} \langle e_1, z_i \rangle \bigr) }  \\
&\qquad \times \one \bigl\{ \, \rho + c^{-1} \langle e_1,z_i \rangle \geq 0, \ i = 1,\dots, k-\ell \,  \bigr\}\, h_t(0,\by, \bz)\, d\bz\, H_\ell (d\rho\, d\by),
\end{align*}
and, accordingly,
\begin{align*}
W_1 (t) &:= \int_0^\infty \int_{(\bbr^d)^{k-1}} e^{-\sum_{i=1}^{k-1} \bigl( \rho + c^{-1} \langle e_1, z_i \rangle \bigr) }  \\
&\qquad \times \one \bigl\{ \, \rho + c^{-1} \langle e_1,z_i \rangle \geq 0, \ i = 1,\dots, k-1 \,  \bigr\}\, h_t(0, \bz)\, d\bz\, H_1 (d\rho), \\
W_k (t) &:= \int_{[0,\infty)\times (\bbr^d)^{k -1}} h_t(0,\by)\, H_k (d\rho\, d\by).
\end{align*}
As we did in Section \ref{s:limit.heavy}, by the decomposition $h_t = h_t^+ - h_t^-$, one can write the process $\BW_\ell$ as the corresponding difference
$\BW_\ell = \BW_\ell^+ - \BW_\ell^-$ for $\ell=1,\dots,k$.

It is easy to compute the covariance function of $\BW_\ell$. We have, for $\ell = 1,\dots,k$ and $t,s\geq0$,
\begin{align}
M_\ell(t,s) &:= \E \bigl\{ W_\ell(t) W_\ell(s) \bigr\} \label{e:cov.comp.light} \\
&=D_\ell \int_0^\infty  \int_{(\bbr^d)^{2k-\ell-1}} \hspace{-20pt}  e^{ -(2k-\ell)\rho - c^{-1} \sum_{i=1}^{2k-\ell-1} \langle e_1,y_i \rangle } \notag  \\
&\qquad \times \one \bigl\{ \rho + c^{-1} \langle e_1, y_i \rangle \geq 0\,, \ i=1,\dots,2k-\ell-1 \bigr\} h_{t,s}^{(\ell)}(0,\by)\, d\by\, d\rho\,, \notag
\end{align}
where 
\begin{equation}  \label{e:def.h.ell}
h_{t,s}^{(\ell)} (0,y_1,\dots,y_{2k-\ell-1}) := h_t(0,y_1,\dots,y_{k-1})\, h_s(0,y_1,\dots,y_{\ell-1,} y_{k}, \dots, y_{2k-\ell-1})\,,
\end{equation}
and, in particular, we set
$$
h_s(0,y_1,\dots,y_{\ell-1}, y_{k}, \dots, y_{2k-\ell-1}) := \begin{cases}
h_s(0,y_{k}, \dots, y_{2k-2})  &  \text{if } \ell = 1\,,  \\
h_s(0,y_{1}, \dots, y_{k-1})  &  \text{if } \ell = k\,.
\end{cases}
$$

It is important to note that if $a(z) \to \infty$ as $z\to\infty$, then $M_\ell$ coincides with $L_\ell$ given in \eqref{e:cov.comp} up to multiplicative factors, i.e.,
$$
M_\ell (t,s) = \bigl( \alpha - d(2k-\ell)^{-1} \bigr) L_\ell(t,s), \ \ t,s \geq0.
$$
This in turn implies that 
$$
\BW_\ell \stackrel{d}{=} \bigl( \alpha - d(2k-\ell)^{-1} \bigr)^{1/2} \BV_\ell,  
$$
in which case, there is nothing to explore here, because the properties of $\BV_\ell$ have  already been studied in Section \ref{s:limit.heavy}. 

In contrast, if $a(z) \to c \in (0,\infty)$ as $z\to\infty$, then $M_\ell$ does not directly relate to $L_\ell$  as above, and, consequently, the process $\BW_\ell$ exhibits properties different to those of $\BV_\ell$. For example, although one may anticipate, as the
analog of the process $\BV_1$, that $\BW_1$ is
a degenerate Gaussian process, this is no longer the case.
\begin{proposition}
Suppose that $a(z) \to c\in(0,\infty)$ as $z\to\infty$. \\
$(i)$ $\BW_1$ is a non-degenerate Gaussian process. \\
$(ii)$ For $\ell = 1,\dots,k$, $\BW_\ell$ is non-self-similar. 
\end{proposition}
\begin{proof}
If $a(z)\to c\in (0,\infty)$ as $z\to \infty$, then $M_1(t,s)$ cannot be decomposed into a function of $t$ and a function of $s$, and therefore, $\BW_1$ is non-degenerate. 

As for $(ii)$, $M_\ell$ does not match $L_\ell$ at all and it loses the scale invariance, meaning that $\BW_\ell$ is non-self-similar. 
\end{proof}

Similarly to Proposition \ref{p:limit.heavy1}, however, the process $\BW_k (= \BW_k^+ - \BW_k^-)$ can be denoted in law as the difference between two time-changed Brownian motions, regardless of whether $a(z) \to \infty$ or $a(z) \to c\in (0,\infty)$ as $z \to \infty$. Furthermore,  the sample paths of $\BW_\ell$ are H\"{o}lder continuous. 
\begin{proposition}
Irrespective of the limit of $a$, the following two results hold. \\
$(i)$ The process $\BW_k^+$ can be represented in law as
$$
\bigl( W_k^+(t), \, t\geq 0 \bigr) \stackrel{d}{=} \biggl( B \Bigl( \int_{[0,\infty) \times (\bbr^d)^{k-1}} \hspace{-10pt}h_t^+(0,\by)\, \mu_k(d\rho\, d\by) \Bigr), \, t \geq 0 \biggr),
$$
where $B$ is the standard Brownian motion. \\
The same statement holds for $\BW_k^-$, by replacing $h_t^+$ with $h_t^-$. 
\vspace{5pt}

\noindent $(ii)$ For $\ell = 1,\dots,k$, and every $T>0$, $\bigl( W_\ell(t), \, 0\leq t \leq T
\bigr)$ has a modification, the sample paths of which are
H\"{o}lder continuous of any order in $[0,1/2)$.
\end{proposition}
\begin{proof}
The proof of $(i)$ is very similar to that in Proposition \ref{p:limit.heavy1}, so we omit it. 
The proof of $(ii)$ is analogous to that in Proposition
\ref{p:self.similar.heavy} $(ii)$; we  have only to show that for
some $C>0$,
\begin{equation*}
\E \Bigl\{ \bigl(W_\ell(t) - W_\ell(s)\bigr)^2 \Bigr\} \leq C(t-s) \ \ \text{for all } 0 \leq s \leq t \leq T\,.
\end{equation*}
Because of the decomposition $\BW_\ell = \BW_\ell^+ - \BW_\ell^-$, it
suffices to prove the above for each $\BW_\ell^+$ and
$\BW_\ell^-$. We  check only the case of $\BW_\ell^+$. We see that
\begin{align*}
\E &\Bigl\{ \bigl(W_\ell^+(t) - W_\ell^+(s) \bigr)^2\Bigr\}  \\
&= \int_{[0,\infty)\times (\bbr^d)^{\ell -1}} \biggl( \int_{(\bbr^d)^{k-\ell}} e^{-\sum_{i=1}^{k-\ell} \bigl( \rho + c^{-1} \langle e_1, z_i \rangle  \bigr)} \one \bigl\{ \, \rho + c^{-1} \langle e_1,z_i \rangle \geq 0, \ i = 1,\dots, k-\ell \,  \bigr\}  \\
&\qquad \qquad \qquad \qquad  \times \bigl( h_t^+(0,\by,\bz) -  h_s^+(0,\by,\bz) \bigr)\, d\bz   \biggr)^2 \mu_\ell (d\rho\, d\by) \\
&\leq D_\ell\, B_\ell^{-1}\, \E \Bigl\{ \bigl( V_\ell^+ (t) - V_\ell^+(s) \bigr)^2 \Bigr\}.
\end{align*}
The rest of the argument is completely the same as Proposition \ref{p:self.similar.heavy} $(ii)$.
\end{proof}

Now, we can state the FCLT of the centered and scaled subgraph counting process
$$
X_n(t) = \tau_n^{-1/2} \bigl( G_n(t) - \E \{ G_n(t) \} \bigr)\,, \ \ \ t \geq 0\,,
$$
where the normalizing sequence $(R_n)$ satisfies
\eqref{e:normalizing.light} and $(\tau_n)$ is defined in
\eqref{e:tau.light}. Interestingly, if $f$ has a subexponential
tail, i.e., $a(z) \to \infty$, then the limiting Gaussian
processes in the theorem below completely coincide (up to
multiplicative constants) with those in Theorem
\ref{t:main.heavy}. When $f$ has an exponential tail, i.e., $a(z)
\to c \in (0,\infty)$, the limiting Gaussian
processes are essentially different from those in Theorem \ref{t:main.heavy}.
The proof of the theorem is presented in Section
\ref{s:proof.light}. For the reader's convenience, we summarize in Tables 1 and 2
the properties of the limiting Gaussian processes in Theorems
\ref{t:main.heavy} and \ref{t:main.light}. These tables indicate that the limiting Gaussian processes  are somewhat special when $f$ has an exponential tail. For example, in this case, the limits always lose self-similarity, regardless of the asymptotics of $nf(R_ne_1)$, whereas, in the regularly varying or the subexponential tail case, the self-similarity is lost only when $nf(R_ne_1)$ converges to a positive and finite constant. Furthermore, when $nf(R_ne_1) \to \infty$, a non-degenerate limit appears only in the exponential tail case. 
\begin{theorem}  \label{t:main.light}
Assume that the density \eqref{e:density.light} satisfies \eqref{e:von.Mises}, \eqref{e:flat}, \eqref{e:poly.upper},
and \eqref{e:subexp.tail}. \\
$(i)$ If $nf(R_ne_1) \to 0$ as $n \to \infty$, then
$$
\bigl( X_n(t), \, t\geq 0 \bigr) \Rightarrow \bigl( W_k(t), \, t\geq 0 \bigr) \ \ \text{in } \mathcal{D}[0,\infty)\,.
$$
\noindent $(ii)$ If $nf(R_ne_1) \to \xi \in (0,\infty)$ as $n \to \infty$, then
$$
\bigl( X_n(t), \, t\geq 0 \bigr) \Rightarrow \left( \sum_{\ell = 1}^k \xi^{2k-\ell} W_\ell(t), \, t\geq 0 \right)  \ \ \text{in } \mathcal{D}[0,\infty)\,.
$$
\noindent $(iii)$ If $nf(R_ne_1) \to \infty$ as $n \to \infty$, then
$$
\bigl( X_n(t), \, t\geq 0 \bigr) \Rightarrow \bigl( W_1(t), \, t\geq 0 \bigr) \ \ \text{in } \mathcal{D}[0,\infty)\,.
$$
\end{theorem}
\begin{table}[htb]
\begin{tabular}{c|ccc} \hline  \\
 & $nf(R_ne_1) \to 0$ & $nf(R_ne_1) \to \xi$ & $nf(R_ne_1) \to \infty$ \\[5pt] \hline \hline
\rule[0pt]{0pt}{18pt}
Regularly varying tail & $d(k-1)/2$ & Non-SS & $d(k-1)$ \\[10pt]
Subexponential tail & $d(k-1)/2$ & Non-SS & $d(k-1)$ \\[10pt]
Exponential tail & Non-SS & Non-SS & Non-SS \\[5pt] \hline \hline
\end{tabular}
\caption{Self-similarity exponents of the limiting Gaussian
processes. Non-SS means that the process is non-self-similar. A
zero limit of $nf(R_ne_1)$ is equivalent to the case in which a
ball $B(0,R_n)$ contains a weak core, and $nf(R_ne_1) \to \infty$
if and only if $B(0,R_n)$ is contained in a weak core. If
$nf(R_ne_1) \to \xi \in (0,\infty)$, then $B(0,R_n)$ agrees with a
weak core (up to multiplicative constants). }
\end{table}
\begin{table}[htb]
\begin{tabular}{c|p{130pt}cp{100pt}} \hline \\
 & $\hspace{30pt} nf(R_ne_1) \to 0$ & $nf(R_ne_1) \to \xi$ & $\hspace{15pt} nf(R_ne_1) \to \infty$ \\[10pt] \hline \hline
\rule[0pt]{0pt}{20pt}
Regularly varying tail & Difference of time-changed Brownian motions & New & Degenerate Gaussian process \\[15pt]
Subexponential tail & Difference of time-changed Brownian motions & New & Degenerate Gaussian process \\[15pt]
 Exponential tail & Difference of time-changed Brownian motions & New & New \\[15pt] \hline \hline
\end{tabular}
% The caption should be at the top of the table. I do not remember how to do this.
\caption{Representation results on the limiting Gaussian
processes. ``New" implies that the limit constitutes a new class
of Gaussian processes.}
\end{table}
%\vspace{10pt}

\section{Graph Connectivity in Weak Core}  \label{s:connect.core}

We start this section by recalling the \textit{weak core}, which
was defined as a centered ball $B \bigl( 0,R_n^{(w)} \bigr)$ such
that $nf \bigl( R_n^{(w)}e_1 \bigr) \to 1$ as $n \to \infty$. In
addition, we need the relevant notion, the \textit{core}, which was defined in Definition \ref{def.core}. Recall that, given a Poisson
point process $\mathcal{P}_n$ on $\bbr^d$, a core is a centered
ball $B(0, R_n)$ such that
\begin{equation}  \label{e:def.core}
B(0, R_n) \subset \bigcup_{X \in \mathcal{P}_n \cap B(0, R_n)} B(X, 1)\,.
\end{equation}

In the following, we seek the largest possible sequence $R_n \to
\infty$ such that the event \eqref{e:def.core} occurs
asymptotically with probability $1$, and subsequently, it is shown
that the largest possible core and a weak core are ``close" in
size. However, the degree of this closeness depends on the tail of
an underlying density $f$, and therefore, we  divide the argument
into two cases.

We first assume that the density $f$ on $\bbr^d$ is spherically
symmetric and has a regularly varying tail, as in
\eqref{e:RV.tail}. For increased clarity, we place an extra
condition that $p(r) := f(re_1)$ is eventually non-increasing in
$r$, that is, $p$ is non-increasing on $(r_0,\infty)$ for some
large $r_0>0$. In this case, the radius of a weak core is,
clearly, given by
\begin{equation}  \label{e:weak.core.heavy}
R_n^{(w)} = \left(\frac{1}{p}\right)^{\leftarrow} (n) := \inf \biggl\{  s: \left( \frac{1}{p} \right)(s) \geq n \biggr\}.
\end{equation}
\begin{proposition}  \label{p:connect.heavy}
Suppose that $p \in RV_{-\alpha}$ for some $\alpha > d$ and $p$ is eventually non-increasing. Define
\begin{equation}  \label{e:max.core.heavy}
R_n^{(c)} = \left(\frac{1}{p}\right)^{\leftarrow} \left( \frac{\delta_1 n}{\log n - \delta_2 \log \log n} \right)
\end{equation}
with $\delta_1 \in \bigl( 0, \alpha/(2^d d^{d/2+1}) \bigr)$ and $\delta_2 \in (0,1)$. If $R_n \leq R_n^{(c)}$, then
\begin{equation}  \label{e:asym.con.heavy}
\P \left( \, B(0, R_n) \subset \bigcup_{X \in \mathcal{P}_n \cap B(0, R_n)} B(X, 1)\, \right) \to 1\,, \ \ \ n \to \infty\,.
\end{equation}
Furthermore, the sequences $(R_n^{(c)})$ in
\eqref{e:max.core.heavy} and $(R_n^{(w)})$ in
\eqref{e:weak.core.heavy} are both regularly varying sequences
with exponent $1/\alpha$, and
\begin{equation}  \label{e:comp.Rc.R0.heavy}
\frac{R_n^{(c)}}{R_n^{(w)}} - \left( \frac{\delta_1}{\log n - \delta_2 \log \log n} \right)^{1/\alpha} \to 0\,, \ \ \ n \to \infty\,.
\end{equation}
\end{proposition}

One can obtain a parallel result when the underlying density has
an exponentially decaying tail, as in \eqref{e:density.light}. We
simplify the situation a bit by assuming
\begin{equation}  \label{e:density.light.sim}
f(x) = C \exp \bigl\{ -\psi \bigl( ||x|| \bigr) \bigr\}, \ \ \ x \in \bbr^d\,,
\end{equation}
where $C$ is a normalizing constant and $\psi: \bbr_+ \to \bbr$ is of
$C^2$-class and satisfies $\psi \in RV_v$ (at infinity) for some
$v>0$ and $\psi^{\prime} >0$. It should be noted that we are
permitting the case $v>1$, implying that, unlike in the previous
section, we do not rule out densities with superexponential tail. Evidently,
the radius of a weak core is given by
\begin{equation}  \label{e:weak.core.light}
R_n^{(w)} = \psiinv(\log n + \log C)\,.
\end{equation}
\begin{proposition}  \label{p:connect.light}
Assume that a probability density $f$ on $\bbr^d$ is given by \eqref{e:density.light.sim}. Define
\begin{equation}  \label{e:max.core.light}
R_n^{(c)} = \psi^{\leftarrow } (\log n - \log \log \log n - \delta_1 - \delta_2)\,,
\end{equation}
where $\delta_1 = d \log 2 - \log v + (1+d/2) \log d - \log C$ and $\delta_2>0$. If $R_n \leq R_n^{(c)}$, then
\begin{equation}  \label{e:asym.con.light}
\P \left( \, B(0, R_n) \subset \bigcup_{X \in \mathcal{P}_n \cap B(0, R_n)} B(X, 1)\, \right) \to 1\,, \ \ \ n \to \infty\,.
\end{equation}
Furthermore, the sequences $(R_n^{(c)})$ in
\eqref{e:max.core.light} and $(R_n^{(w)})$ in
\eqref{e:weak.core.light} are close in size in the sense of
\begin{equation}  \label{e:comp.Rc.R0.light}
\frac{R_n^{(c)}}{R_n^{(w)}} - \left( 1-\frac{\log \log \log n + \delta_1 + \delta_2 + \log C}{\log Cn} \right)^{1/v} \to 0\,, \ \ \ n \to \infty\,.
\end{equation}
\end{proposition}
The following result is needed as preparation for the proof of
these propositions. The proof may be obtained by slightly
modifying the proof of Theorem 2.1 in
\cite{adler:bobrowski:weinberger:2014}, but so that this paper is
self-contained, we repeat the argument.
\begin{lemma}  \label{l:Qlemma}
Given a spherically symmetric density $f$ on $\bbr^d$, suppose
that $p(r) = f(re_1)$ is eventually non-increasing. Let $g = 1/(2
d^{1/2})$. Suppose, in addition, that there exists a sequence $R_n
\nearrow \infty$ such that $d \log R_n - g^d n f(R_ne_1) \to -\infty$
as $n \to \infty$. Then,
\begin{equation}  \label{e:only.lower}
\P \left( \, B(0, R_n) \subset \bigcup_{X \in \mathcal{P}_n \cap B(0, R_n)} B(X, 1)\, \right) \to 1\,, \ \ \ n \to \infty\,.
\end{equation}
\end{lemma}
\begin{proof}
For $\rho>0$, let $\mathcal{Q}(\rho)$ be a collection of cubes
with grid $g$ that are contained in $B(0, \rho)$. Then,
$$
\bigl\{ \, Q \cap \mathcal{P}_n \neq \emptyset \ \, \text{for all } Q \in \mathcal{Q}(\rho)\, \bigr\} \subset \Bigl\{ \, B(0, \rho) \subset \bigcup_{X \in \mathcal{P}_n \cap B(0, \rho)} B(X, 1)\, \Bigr\}
$$
for all $\rho>0$ and $n \geq 1$. It now suffices to show that
$$
\P \bigl( \, Q \cap \mathcal{P}_n = \emptyset \ \, \text{for some } Q \in \mathcal{Q}(R_n)\, \bigr) \to 0\,, \ \ \ n \to \infty\,.
$$
This probability is estimated from above by
$$
\sum_{Q \in \mathcal{Q}(R_n)} \P(Q \cap \mathcal{P}_n = \emptyset) = \sum_{Q \in \mathcal{Q}(R_n)} \exp \Bigl\{ -n \int_Q f(x) dx \Bigr\}
$$
$$
\leq \sum_{Q \in \mathcal{Q}(R_n)} \exp \bigl\{ -n g^d f(R_n e_1) \bigr\} \leq g^{-d} R_n^d \exp \bigl\{ -g^d n f(R_ne_1) \bigr\}\,.
$$
At the first inequality, we used the fact that $p$ is eventually
non-increasing. Clearly, the rightmost term vanishes as $n \to
\infty$.
\end{proof}
\begin{proof}(\textit{proof of Proposition \ref{p:connect.heavy}})
Observe that the assumption $p \in RV_{-\alpha}$ implies
$(1/p)^{\leftarrow} \in RV_{1/\alpha}$, e.g., Proposition 2.6 (v)
in \cite{resnick:2007}. Thus, \eqref{e:comp.Rc.R0.heavy} readily
follows from the uniform convergence of regularly varying
functions; see Proposition 2.4 in \cite{resnick:2007}. By Lemma
\ref{l:Qlemma}, it suffices to verify that $d \log R_n^{(c)} - g^d
nf \bigl( R_n^{(c)}e_1 \bigr) \to -\infty$ as $n \to \infty$.
Since $0 < \delta_2 < 1$, we have
\begin{align*}
d \log R_n^{(c)} &\leq d  \left[ \log \left(\frac{1}{p}\right)^{\leftarrow}  \left( \frac{\delta_1 n}{\log n - \delta_2 \log \log n} \right) \right] \left[ \log \left( \frac{\delta_1 n}{\log n - \delta_2 \log n \log n} \right) \right]^{-1}  \\
&\times (\log \delta_1 + \log n - \delta_2 \log \log n)\,,
\end{align*}
and $g^d nf \bigl( R_n^{(c)}e_1 \bigr) = g^d \delta_1^{-1} (\log n
- \delta_2 \log \log n)$. Using Proposition 2.6 (i) in
\cite{resnick:2007},
\begin{align*}
&d \left[ \log \left(\frac{1}{p}\right)^{\leftarrow}  \left( \frac{\delta_1 n}{\log n - \delta_2 \log \log n} \right) \right] \left[ \log \left( \frac{\delta_1 n}{\log n - \delta_2 \log n \log n} \right) \right]^{-1}  -  g^d \delta_1^{-1} \\
&\to d \alpha^{-1} -  g^d \delta_1^{-1} < 0\,, \ \ \ n \to \infty\,.
\end{align*}
At the last inequality, we applied the constraint in $\delta_1$.
Therefore, we have $d \log R_n^{(c)} - g^d nf \bigl( R_n^{(c)}e_1
\bigr) \to -\infty$, $n \to \infty$, as requested.
\end{proof}
\begin{proof}(\textit{proof of Proposition \ref{p:connect.light}})
Since $\psiinv \in RV_{1/v}$, it is easy to show
\eqref{e:comp.Rc.R0.light}, and therefore, we  prove only that $d
\log R_n^{(c)} - g^d nf \bigl( R_n^{(c)}e_1 \bigr) \to -\infty$ as
$n \to \infty$. We see that
$$
d \log R_n^{(c)} \leq d \log \psiinv (\log n) \sim d v^{-1} \log \log n\,, \ \ \ n \to \infty\,,
$$
and that $g^d n f \bigl( R_n^{(c)}e_1 \bigr)= g^d C e^{\delta_1 +
\delta_2} \log \log n$. By virtue of the constraints in $\delta_1$
and $\delta_2$, we have $dv^{-1} - g^d C e^{\delta_1 + \delta_2} <
0$; thus, the claim is proved.
\end{proof}
\begin{remark}
The proof of Lemma \ref{l:Qlemma} merely estimated the probability
in \eqref{e:only.lower} from below. Therefore, it seems to be
possible that in the propositions above, \eqref{e:asym.con.heavy}
and \eqref{e:asym.con.light} may hold for the sequence $R_n \nearrow
\infty$ growing more quickly than $R_n^{(c)}$ but more slowly than
$R_n^{(w)}$, i.e., $R_n^{(c)} \leq R_n \leq R_n^{(w)}$; it is
unknown, however, to what extent we can make $R_n$ closer to
$R_n^{(w)}$.
\end{remark}
%\vspace{10pt}

\section{Proof of Main Results}  \label{s:proof.main}

This section presents the proof of the main results of this paper.
The proof is, however, rather long, and therefore, it is divided
into several parts. All the supplemental ingredients necessary are
collected in the Appendix, most of which are cited from
\cite{penrose:2003}.

Let Ann$(K,L)$ be an annulus of inner radius $K$ and outer radius
$L$. For $x_1,\dots,x_k \in \bbr^d$, define $\text{Max}(x_1,\dots,
x_k)$ as the function selecting an element with the largest
distance from the origin. That is, $\text{Max}(x_1,\dots, x_k) =
x_i$ if $||x_i|| = \max_{1 \leq j \leq k}||x_j||$. If multiple
$x_j$'s achieve the maximum, we choose an element with the
smallest subscript.

In the following, $\Y, \Y^{\prime}, \Y_i$, etc. always represent a
finite collection of $d$-dimensional real vectors. We use the
following shorthand notations. That is, for $\bx = (x_1,\dots,x_m)
\in (\bbr^d)^m$, $x \in \bbr^d$, and $\by = (y_1,\dots,y_{m-1})
\in (\bbr^d)^{m-1}$,
\begin{align*}
f(\bx) &:= f(x_1)\cdots f(x_m)\,, \\
f(x+\by) &:= f(x + y_1)\cdots f(x + y_{m-1})\,, \\
h(0,\by) &:= h(0,y_1,\dots,y_{m-1}) \ \ \text{etc}.
\end{align*}

Regarding the indicator $h_t :(\bbr^d)^k \to \{ 0,1 \}$ given in \eqref{e:geo.graph.dyna}, the following notations are used to save
space.
\begin{align}
h_{t,s} (\bx) &:= h_t(\bx) - h_s (\bx)\,, \ \ 0 \leq s \leq t\,, \  \bx \in (\bbr^d)^k,  \label{e:def.hts} \\
h_{t,s}^{\pm} (\bx) &:= h_t^{\pm}(\bx) - h_s^{\pm} (\bx)\,, \ \ 0 \leq s \leq t\,, \  \bx \in (\bbr^d)^k,  \notag \\
h_{n,t,s} (\bx) &:= h_{t,s} (\bx)\, \one \bigl\{ m(\bx) \geq R_n \bigr\}\,, \ \ 0 \leq s \leq t\,, \  \bx \in (\bbr^d)^k,  \label{e:def.hnts}
\end{align}
and for $\ell \in \{ 0,\dots,k\}$,
\begin{equation}  \label{e:def.htsl}
h_{t,s}^{(\ell)}(\bx) := h_t(x_1,\dots,x_k)\, h_s(x_1,\dots,x_\ell, x_{k+1}, \dots, x_{2k-\ell})\,,  \ \ t,s \geq 0\,, \ \bx \in (\bbr^d)^{2k-\ell}.
\end{equation}
In particular, we set
$$
h_s(x_1,\dots,x_\ell, x_{k+1}, \dots, x_{2k-\ell}) := \begin{cases}
h_s(x_{k+1}, \dots, x_{2k})  &  \text{if } \ell = 0\,,  \\
h_s(x_{1}, \dots, x_{k})  &  \text{if } \ell = k\,.
\end{cases}
$$
In Section \ref{s:proof.heavy}, we use, for $1 \leq K < L \leq \infty$, $n \in \bbn_+$ and $t \geq 0$,
\begin{align*}
h_{n,t,K,L}(\bx) &:= h_t (\bx)\, \one \bigl\{ m(\bx) \geq R_n\,, \ \text{Max}(\bx) \in \text{Ann}(KR_n,LR_n) \}\,, \\
h_{n,t,K,L}^{\pm}(\bx) &:= h_t ^{\pm}(\bx)\, \one \bigl\{ m(\bx) \geq R_n\,, \ \text{Max}(\bx) \in \text{Ann}(KR_n,LR_n) \}\,.
\end{align*}
The same notations are retained for Section \ref{s:proof.light} to
represent, for $0 \leq K < L \leq \infty$, $n \in \bbn_+$ and $t \geq
0$,
\begin{align*}
h_{n,t,K,L}(\bx) &:= h_t (\bx)\, \one \bigl\{ m(\bx) \geq R_n\,, \ a(R_n)^{-1}\bigl(\text{Max}(\bx) - R_n \bigr) \in [K,L) \bigr\}\,,   \\
h_{n,t,K,L}^{\pm}(\bx) &:= h_t^{\pm} (\bx)\, \one \bigl\{ m(\bx) \geq R_n\,, \ a(R_n)^{-1}\bigl(\text{Max}(\bx) - R_n \bigr) \in [K,L) \bigr\}\,.
\end{align*}

Finally, $C^*$ denotes a generic positive constant, which may
change between lines and does not depend on $n$.

In the following, we divide the argument into two subsections.
Section \ref{s:proof.heavy} treats the case in which the
underlying density has a regularly varying tail; our goal is to
prove Theorem \ref{t:main.heavy}. Subsequently Section
\ref{s:proof.light} provides the proof of Theorem
\ref{t:main.light}, where the density is assumed to have an
exponentially decaying tail. Before the specific subsections,
however, we show some preliminary results, which are commonly used
in both subsections for the tightness proof.

\begin{lemma}  \label{l:used.for.tightness}
Let $h_t: (\bbr^d)^k \to \{  0,1\} $ be an indicator given in  \eqref{e:geo.graph.dyna}. Fix $T > 0$. Then, we have for $\ell \in \{ 1,\dots,k\}$,
\begin{align}  
\int_{(\bbr^d)^{\ell-1}} \hspace{-10pt} d\by \int_{(\bbr^d)^{k-\ell}} \hspace{-10pt} d\bz_2 \int_{(\bbr^d)^{k-\ell}} \hspace{-10pt} d\bz_1 h_{t,s}^{+}(0,\by, \bz_1)\, h_{s,r}^{+}(0,\by, \bz_2)
&\leq C^* (t - s) (s - r)\,, \label{e:triple.int} \\
\int_{(\bbr^d)^{\ell-1}} \hspace{-10pt} d\by \int_{(\bbr^d)^{k-\ell}} \hspace{-10pt} d\bz_2 \int_{(\bbr^d)^{k-\ell}} \hspace{-10pt} d\bz_1 h_{t,s}^{-}(0,\by, \bz_1)\, h_{s,r}^{-}(0,\by, \bz_2)
&\leq C^* (t - s) (s - r) \notag 
\end{align}
for all $0 \leq r \leq s \leq t \leq T$.
\end{lemma}
\begin{proof}
We only prove the first inequality. If $\ell = 1$ or $\ell = k$, the claim is trivial, and therefore,
we can take $2 \leq \ell \leq k-1$. It follows from
\eqref{e:close.enough.decomp.dyna} that the integral in \eqref{e:triple.int}
is not altered if the integral domain is restricted to $\bigl(
B(0,kT) \bigr)^{\ell-1} \times \bigl( B(0,kT) \bigr)^{k-\ell}
\times \bigl( B(0,kT) \bigr)^{k-\ell}$. With $\lambda$ being the
Lebesgue measure on $(\bbr^d)^{k-\ell}$, we see that for every
$\by \in (\bbr^d)^{\ell -1}$,
\begin{align}
\int_{ \bigl( B(0,kT) \bigr)^{k-\ell}} h_{t,s}^+(0,\by, \bz)  d\bz &= \lambda \bigl\{ \bz \in \bigl( B(0,kT) \bigr)^{k-\ell}:  h_t^+(0,\by,\bz) = 1\,, \ h_s^+(0,\by,\bz) = 0  \bigr\}  \label{e:seq.volume}  \\
&\leq  \lambda \bigl\{ \bz \in \bigl( B(0,kT) \bigr)^{k-\ell}:  s < ||z_i - z_j|| \leq t \ \text{for some } i \neq j \bigr\}  \notag \\
&\quad + \lambda \bigl\{ \bz \in \bigl( B(0,kT) \bigr)^{k-\ell}:  s < ||z_i - y_j|| \leq t \ \text{for some } i, j \bigr\}  \notag \\
&\quad + \lambda \bigl\{ \bz \in \bigl( B(0,kT) \bigr)^{k-\ell}:  s < ||z_i|| \leq t \ \text{for some } i \bigr\} \notag \\
&\quad + \lambda \bigl\{ \bz \in \bigl( B(0,kT) \bigr)^{k-\ell}:  s < ||y_i - y_j|| \leq t \ \text{for some } i \neq j \bigr\} \notag \\
&\quad + \lambda \bigl\{ \bz \in \bigl( B(0,kT) \bigr)^{k-\ell}:  s < ||y_i|| \leq t \ \text{for some } i \bigr\}.  \notag
\end{align}
Observe that for $i \neq j$,
$$
\lambda \bigl\{ \bz \in \bigl( B(0,kT) \bigr)^{k-\ell}:  s < ||z_i - z_j|| \leq t  \bigr\} \leq (kT)^{d(k-\ell-1)} (\omega_d)^{k-\ell} (t^d-s^d)\,,
$$
where $\omega_d$ is the volume of the $d$-dimensional unit ball.
Since the second and the third terms on the rightmost term in
\eqref{e:seq.volume} have the same upper bound, we ultimately
obtain
\begin{align*}
\int_{ \bigl( B(0,kT) \bigr)^{k-\ell}} &h_{t,s}^{+}(0,\by, \bz)  d\bz \\
&\leq C^* \biggl(  t^d-s^d + \sum_{i,j=1, \ i \neq j}^{\ell - 1} \one \bigl\{ s < ||y_i - y_j|| \leq t \bigr\} + \sum_{i=1}^{\ell-1}  \one \bigl\{ s < ||y_i|| \leq t \bigr\}\biggr).
\end{align*}
Therefore, the integral in \eqref{e:triple.int} is bounded above by
\begin{align*}
C^* \int_{\bigl( B(0,kT) \bigr)^{\ell-1}} &\biggl(  t^d-s^d + \sum_{i,j=1, \ i \neq j}^{\ell - 1} \one \bigl\{ s < ||y_i - y_j|| \leq t \bigr\} + \sum_{i=1}^{\ell-1}  \one \bigl\{ s < ||y_i|| \leq t \bigr\}\biggr) \\
&\times \biggl(  s^d-r^d + \sum_{i,j=1, \ i \neq j}^{\ell - 1} \one \bigl\{ r < ||y_i - y_j|| \leq s \bigr\} + \sum_{i=1}^{\ell-1}  \one \bigl\{ r < ||y_i|| \leq s \bigr\}\biggr) d\by\,.
\end{align*}
An elementary calculation shows that for all $i,j, i^{\prime},
j^{\prime} \in \{ 1,\dots,\ell-1 \}$ with $i>j$ and
$i^{\prime}>j^{\prime}$,
\begin{align*}
&\int_{\bigl( B(0,kT) \bigr)^{\ell-1}} \one \bigl\{ s < ||y_i - y_j|| \leq t \bigr\}\, \one \bigl\{ r < ||y_{i^{\prime}} - y_{j^{\prime}}|| \leq s \bigr\} d\by \\
&\quad \leq C^* (t^d-s^d)(s^d-r^d) \leq C^* (t-s)(s-r)
\end{align*}
In particular, if $i=i^{\prime}$ and $j=j^{\prime}$, the integral
is identically zero. Applying the same manipulation to the
integral of other cross-terms, we can conclude the claim of the
lemma.
\end{proof}

\subsection{Regularly Varying Tail Case}  \label{s:proof.heavy}
Under the setup of Theorem \ref{t:main.heavy}, we first define the
subgraph counting process with restricted domain. For $1 \leq K <
L \leq \infty$, $n \in \bbn_+$, and $t \geq 0$, let
\begin{align*}
G_{n,K,L}(t) &= \sum_{\Y \subset \mathcal{P}_n} h_t (\Y)\, \one \bigl\{ m(\Y) \geq R_n\,, \ \text{Max}(\Y) \in \text{Ann}(KR_n,LR_n) \} \\
&:= \sum_{\Y \subset \mathcal{P}_n} h_{n,t,K,L} (\Y)\,,
\end{align*}
and
\begin{align*}
G_{n,K,L}^{\pm}(t) &= \sum_{\Y \subset \mathcal{P}_n} h_t^{\pm} (\Y)\, \one \bigl\{ m(\Y) \geq R_n\,, \  \text{Max}(\Y) \in \text{Ann}(KR_n,LR_n) \} \\
&:= \sum_{\Y \subset \mathcal{P}_n} h_{n,t,K,L}^{\pm} (\Y)\,,
\end{align*}
where $(R_n)$ satisfies \eqref{e:normalizing.heavy}. For the
special case $K=1$ and $L=\infty$, we simply denote $G_n(t) =
G_{n,1,\infty}(t)$ and $G_n^{\pm}(t) = G_{n,1,\infty}^{\pm}(t)$.
The subgraph counting processes, centered and scaled, for which we
prove the FCLT, are given by
\begin{align}
X_n(t) &= \tau_n^{-1/2} \Bigl( G_n(t) - \E\bigl\{ G_n(t) \bigr\} \Bigr)\,,  \notag \\
X_n^{\pm}(t) &= \tau_n^{-1/2} \Bigl( G_n^{\pm}(t) - \E\bigl\{ G_n^{\pm}(t) \bigr\} \Bigr)\,,  \label{e:def.Xnpm.heavy}
\end{align}
where $(\tau_n)$ is determined by \eqref{e:tau.heavy} according to
which regime is considered. The first proposition below computes
the covariances of $\bigl( G_{n,K,L}(t) \bigr)$.
\begin{proposition}  \label{p:cova.heavy}
Assume the conditions of Theorem \ref{t:main.heavy}. Let $1 \leq K < L \leq \infty$. \\
$(i)$ If $nf(R_ne_1) \to 0$ as $n \to \infty$, then
$$
\tau_n^{-1} \text{Cov} \bigl( G_{n,K,L}(t), G_{n,K,L}(s) \bigr) \to (K^{d-\alpha k} - L^{d - \alpha k}) L_k(t,s)\,, \ \ \ n\to \infty\,.
$$
$(ii)$ If $nf(R_ne_1) \to \xi \in (0,\infty)$ as $n \to \infty$, then
$$
\tau_n^{-1} \text{Cov} \bigl( G_{n,K,L}(t), G_{n,K,L}(s) \bigr) \to \sum_{\ell=1}^k (K^{d-\alpha (2k-\ell)} - L^{d - \alpha (2k-\ell)}) \xi^{2k-\ell}L_\ell(t,s)\,, \ \ \ n\to \infty\,.
$$
$(iii)$ If $nf(R_ne_1) \to \infty$ as $n \to \infty$, then
$$
\tau_n^{-1} \text{Cov} \bigl( G_{n,K,L}(t), G_{n,K,L}(s) \bigr) \to (K^{d-\alpha (2k-1)} - L^{d - \alpha (2k-1)}) L_1(t,s)\,, \ \ \ n\to \infty\,.
$$
\end{proposition}
\begin{proof}
We start by writing
\begin{align*}
\E \bigl\{ &G_{n,K,L}(t)\, G_{n,K,L}(s) \bigr\} \\
&= \sum_{\ell=0}^k \E \Bigl\{ \, \sum_{\Y_1 \subset \mathcal{P}_n} \sum_{\Y_2 \subset \mathcal{P}_n} h_{n,t, K, L} (\Y_1)\, h_{n,s, K, L} (\Y_2)\, \one \bigl\{ \, |\Y_1 \cap \Y_2| = \ell\, \bigr\} \Bigr\} \\
&:= \sum_{\ell=0}^k \E \{I_\ell\}\,.
\end{align*}
For $\ell=0$, applying Palm theory (see the Appendix) twice,
\begin{align*}
\E \{I_0\} &= \frac{n^{2k}}{(k!)^2}\, \E \bigl\{ h_{n,t, K, L} (X_1,\dots,X_k)\, h_{n,s, K, L} (X_{k+1},\dots,X_{2k})\bigr\} \\
&= \E \bigl\{G_{n,K,L}(t) \bigr\} \E \bigl\{ G_{n,K,L}(s)\bigr\}\,.
\end{align*}
Therefore, the multiple applications of Palm theory yield
\begin{align*}
\text{Cov} &\bigl( G_{n,K,L}(t)\,, G_{n,K,L}(s) \bigr) = \sum_{\ell=1}^k \E\{ I_\ell\} \\
&= \sum_{\ell=1}^k \frac{n^{2k-\ell}}{\ell ! \bigl( (k-\ell)! \bigr)^2}\, \E \Bigl\{h_{n,t, K, L} (\Y_1)\, h_{n,s, K, L} (\Y_2)\, \one \bigl\{ \, |\Y_1 \cap \Y_2| = \ell\, \bigr\} \Bigr\}.
\end{align*}
Define for $\ell \in \{1,\dots,k  \}$, 
\begin{multline*}
C_n^{(\ell)} (K,L) := \bigl\{ \bx \in (\bbr^d)^{2k-\ell}: \text{Max}(x_1,\dots, x_k) \in \text{Ann}(KR_n,LR_n)\,, \\
 \text{Max}(x_1,\dots, x_\ell, x_{k+1}, \dots, x_{2k-\ell}) \in \text{Ann}(KR_n,LR_n)\bigr\}\,.
\end{multline*}
By the change of variables $\bx \rightarrow (x,x+\by)$ with $\bx
\in (\bbr^d)^{2k-\ell}$, $x \in \bbr^d$, $\by \in
(\bbr^d)^{2k-\ell-1}$, together with invariance
\eqref{e:location.inv}, while recalling notation
\eqref{e:def.htsl},
\begin{align*}
\E \Bigl\{ &h_{n,t,K,L}(\Y_1)\, h_{n,s,K,L}(\Y_2)\, \one \bigl\{ \, |\Y_1 \cap \Y_2| = \ell\, \bigr\} \Bigr\} \\
&= \int_{(\bbr^d)^{2k-\ell}} f(\bx)\, \one \bigl\{ m(\bx) \geq R_n \bigr\}\,  h_{t,s}^{(\ell)} (\bx)\,\one \bigl\{ \bx \in C_n^{(\ell)}(K,L) \bigr\} d\bx  \\
&= \int_{\bbr^d} \int_{(\bbr^d)^{2k-\ell-1}} f(x)\, f(x+\by)\, \one \bigl\{ m(x, x + \by) \geq R_n \bigr\}\, h_{t,s}^{(\ell)} (0, \by)\,  \\
&\quad \times \one \bigl\{ (x, x+\by) \in C_n^{(\ell)}(K,L) \bigr\} d\by dx\,.
\end{align*}
The polar coordinate transform $x \rightarrow (r,\theta)$ and an
additional change of variable $\rho \rightarrow r/R_n$ yield
\begin{align}
\E \Bigl\{ &h_{n,t,K,L}(\Y_1)\, h_{n,s,K,L}(\Y_2)\, \one \bigl\{ \, |\Y_1 \cap \Y_2| = \ell\, \bigr\} \Bigr\}  \label{e:after.polar} \\
&= R_n^d f(R_ne_1)^{2k-\ell} \int_{S_{d-1}} \hspace{-7pt}J(\theta)d\theta \int_1^{\infty} d\rho \int_{(\bbr^d)^{2k-\ell-1}} \hspace{-5pt} d\by\, \rho^{d-1} \frac{f(R_n\rho e_1)}{f(R_ne_1)} \notag \\
&\quad \times \prod_{i=1}^{2k-\ell-1} \frac{f \bigl( ||R_n \rho \theta + y_i|| e_1 \bigr)}{f(R_ne_1)}\, \one \bigl\{ ||\rho \theta + y_i / R_n|| \geq 1 \bigr\}\, h_{t,s}^{(\ell)}(0,\by) \notag \\
&\quad \times \one \bigl\{ (R_n \rho \theta, R_n \rho \theta + \by) \in C_n^{(\ell)}(K,L) \bigr\} \,, \notag
\end{align}
where $S_{d-1}$ denotes the $(d-1)$-dimensional unit sphere in
$\bbr^d$ and $J(\theta)$ is the usual Jacobian
$$
J(\theta) = \sin^{k-2} (\theta_1) \sin^{k-3} (\theta_2) \cdots \sin (\theta_{k-2})\,.
$$
Note that by the regular variation of $f$
(with exponent $-\alpha$), for every $\rho > 1$, $\theta \in
S_{d-1}$, and $y_i$'s,
\begin{equation}  \label{e:RV.result}
\frac{f(R_n\rho e_1)}{f(R_ne_1)} \to \rho^{-\alpha}\,, \quad \ \  \prod_{i=1}^{2k-\ell-1} \frac{f \bigl( ||R_n \rho \theta + y_i|| e_1 \bigr)}{f(R_ne_1)} \to \rho^{-\alpha(2k-\ell-1)}\,, \ \ \ n \to \infty
\end{equation}
and, furthermore,
\begin{equation}  \label{e:conv.indicator}
\one \bigl\{ (R_n \rho \theta, R_n \rho \theta + \by) \in C_n^{(\ell)}(K,L) \bigr\} \to \one \{ K \leq \rho \leq L \}\,, \ \ \ n \to \infty\,.
\end{equation}
Substituting  \eqref{e:RV.result} and \eqref{e:conv.indicator}
back into \eqref{e:after.polar}, while supposing
temporarily that the dominated convergence theorem is applicable,
we may conclude that
\begin{align}
\text{Cov} &\bigl( G_{n,K,L}(t)\,, G_{n,K,L}(s) \bigr) \label{e:cov.total} \\
&\sim \sum_{\ell=1}^k n^{2k-\ell} R_n^d f(R_ne_1)^{2k-\ell} \bigl( K^{d-\alpha(2k-\ell)} - L^{d-\alpha(2k-\ell)} \bigr) L_\ell(t,s)\,, \ \ \ n \to \infty\,. \notag
\end{align}
Observe that the limit value of $nf(R_ne_1)$ completely determines
which term on the right hand side of \eqref{e:cov.total} is
dominant. If $nf(R_ne_1) \to 0$, then the $k$th term,i.e., $\ell =
k$, in the sum grows fastest, while the first term, i.e., $\ell =
1$, grows fastest when $nf(R_ne_1) \to \infty$. Moreover, if
$nf(R_ne_1) \to \xi \in (0,\infty)$, then all the terms in the sum
grow at the same rate. This concludes the claim of the
proposition.

It now remains to establish an integrable upper bound for the
application of the dominated convergence theorem. First, condition \eqref{e:close.enough} provides
$$
h_{t,s}^{(\ell)}(0,\by) \leq \one \bigl\{ ||y_i|| \leq k(t+s)\,, \ i=1,\dots,2k-\ell-1 \bigr\}\,.
$$
Next, appealing to
Potter's bound ,e.g., Proposition 2.6 $(ii)$ in
\cite{resnick:2007}, for every $\xi \in (0,\alpha-d)$ and
sufficiently large $n$,
$$
\frac{f(R_{n} \rho e_1)}{f(R_{n} e_1)}\, \one \{ \rho \geq 1 \} \leq (1+\xi) \, \rho^{-\alpha +\xi}\, \one \{ \rho \geq 1 \}
$$
and
\begin{align*}
\prod_{i=1}^{2k-\ell-1} \frac{f(||R_n \rho \theta + y_i || e_1)}{f(R_n e_1)}\, \one \bigl\{|| \rho \theta + y_i/R_{n} || \geq 1 \bigr\} \leq (1+\xi)^{2k-\ell-1}\,. \\
\end{align*}
Since $\int_{1}^{\infty} \rho^{d-1-\alpha+\xi} d\rho < \infty$, we
are allowed to apply the dominated convergence theorem.
\end{proof}
The next proposition proves the weak convergence of Theorem
\ref{t:main.heavy} in a finite-dimensional sense.
\begin{proposition}   \label{p:fidi.heavy}
Assume the conditions of Theorem \ref{t:main.heavy}. Then, weak
convergences $(i) - (iii)$ in the theorem hold in a
finite-dimensional sense.
Furthermore, let $\BX_n^{\pm}$ be the processes defined in \eqref{e:def.Xnpm.heavy}. Then, the following results
also hold in a finite-dimensional sense. \\
$(i)$ If $nf(R_ne_1) \to 0$ as $n \to \infty$, then
\begin{equation}  \label{e:bivariate1}
(\BX_n^+, \BX_n^-) \Rightarrow (\BV_k^+, \BV_k^-)\,.
\end{equation}
$(ii)$ If $nf(R_ne_1) \to \xi \in (0,\infty)$ as $n \to \infty$, then
\begin{equation}  \label{e:bivariate2}
(\BX_n^+, \BX_n^-) \Rightarrow \left(\sum_{\ell=1}^k \xi^{2k-\ell} \BV_\ell^+, \, \sum_{\ell=1}^k \xi^{2k-\ell} \BV_\ell^- \right)\,.
\end{equation}
$(iii)$ If $nf(R_ne_1) \to \infty$ as $n \to \infty$, then
\begin{equation}  \label{e:bivariate3}
(\BX_n^+, \BX_n^-) \Rightarrow (\BV_1^+, \BV_1^-)\,.
\end{equation}
The limiting Gaussian processes $(\BV_\ell^+, \BV_\ell^-)$, $\ell
= 1,\dots,k$ are all formulated in Section \ref{s:limit.heavy}.
\end{proposition}
\begin{proof}
The proofs of \eqref{e:bivariate1}, \eqref{e:bivariate2}, and
\eqref{e:bivariate3} are a bit more technical, but are very
similar to the corresponding results in Theorem
\ref{t:main.heavy}; therefore, we  check only finite-dimensional
weak convergences in Theorem \ref{t:main.heavy}. The argument here
is closely related to that in Theorem 3.9 of \cite{penrose:2003},
for which we rely on the so-called Cram\'er-Wold device. For $0
\leq t_1 < \dots < t_m < \infty$, $a_1,\dots, a_m \in \bbr$ and $m
\geq 1$, define $S_n := \sum_{j=1}^m a_j G_n (t_j)$. For $K>1$,
$S_n$ can be further decomposed into two parts:
\begin{align*}
S_n &= \sum_{j=1}^m a_j G_{n,1,K}(t_j) + \sum_{j=1}^m a_j G_{n,K,\infty}(t_j)  \\
&:= T_n^{(K)} + U_n^{(K)}\,.
\end{align*}
We define a constant $\gamma_K$ as follows in accordance with the limit of $nf(R_ne_1)$.
$$
\gamma_K := \begin{cases} \sum_{i=1}^m \sum_{j=1}^m a_i a_j (1-K^{d-\alpha k}) L_k(t_i,t_j) & \text{if } nf(R_ne_1) \to 0\,, \\
\sum_{i=1}^m \sum_{j=1}^m a_i a_j \sum_{\ell=1}^k  (1-K^{d-\alpha(2k-\ell)}) \xi^{2k-\ell} L_\ell(t_i,t_j) & \text{if } nf(R_ne_1) \to \xi \in (0,\infty)\,, \\
\sum_{i=1}^m \sum_{j=1}^m a_i a_j (1-K^{d-\alpha (2k-1)}) L_1(t_i,t_j) & \text{if } nf(R_ne_1) \to \infty\,.
\end{cases}
$$
Moreover, $\gamma := \lim_{K \to \infty} \gamma_K$.
It is then elementary to check that, regardless of the regime we consider,
$$
\tau_n^{-1} \text{Var} \{ T_n^{(K)} \} \to \gamma_K\,, \ \ \ \tau_n^{-1} \text{Var} \{ U_n^{(K)} \} \to \gamma - \gamma_K \ \ \ \text{as } n \to \infty\,.
$$
For the completion of the proof, we ultimately need to show that
$$
\tau_n^{-1/2} \bigl( S_n - \E \{ S_n \} \bigr) \Rightarrow N(0,\gamma)\,.
$$
By the standard approximation argument given on p. 64 of
\cite{penrose:2003}, it suffices to show that
\begin{equation}  \label{e:int.CLT}
\tau_n^{-1/2} \Bigl( T_n^{(K)} - \E \bigl\{ T_n^{(K)} \bigr\} \Bigr) \Rightarrow \text{N}(0, \gamma_K) \ \ \text{for every } K > 1\,;
\end{equation}
equivalently,
\begin{equation}  \label{e:int.normalized.CLT}
\frac{T_n^{(K)} - \E \bigl\{ T_n^{(K)} \bigr\}}{\sqrt{\text{Var}  \bigl\{ T_n^{(K)} \bigr\}}} \Rightarrow N(0,1) \ \ \text{for every } K > 1\,.
\end{equation}

Let $(Q_\ell: \ell \in \bbn)$ be a collection of unit cubes covering $\bbr^d$. Define
$$
V_n := \bigl\{ \ell \in \bbn: Q_\ell \cap \text{Ann}(R_n,KR_n) \neq \emptyset \bigr\}\,,
$$
where we have that $|V_n| \leq C^* R_n^d$.  \\
Then, $T_n^{(K)}$ can be partitioned as follows.
\begin{align*}
T_n^{(K)} &= \sum_{\ell \in V_n} \sum_{j=1}^m a_j \sum_{\Y \subset \mathcal P_n} h_{t_j} (\Y)\, \one \bigl\{ m(\Y) \geq R_n, \ \text{Max}(\Y) \in \text{Ann}(R_n,KR_n) \cap Q_\ell \bigr\} \\
&:= \sum_{\ell \in V_n} \eta_{\ell,n}\,.
\end{align*}

For $i,j \in V_n$, we put an edge between $i$ and $j$ (write $i
\sim j$) if $i \neq j$ and the distance between $Q_i$ and $Q_j$
are less than $2kt_m$. Then, $(V_n, \sim)$ gives a
\textit{dependency graph} with respect to $(\eta_{\ell, n}, \,
\ell \in V_n)$; that is, for any two disjoint subsets $I_1$, $I_2$
of $V_n$ with no edges connecting $I_1$ and $I_2$, $(\eta_{\ell,
n}, \, \ell \in I_1)$ is independent of $(\eta_{\ell, n}, \, \ell
\in I_2)$. Notice that the maximum degree of $(V_n, \sim)$ is at
most finite.

According to Stein's method for normal approximation (see Theorem 2.4 in \cite{penrose:2003}), \eqref{e:int.normalized.CLT} immediately follows if we
can show that for $p=3,4$,
\begin{equation}  \label{e:normal.approx}
R_n^d\, \max_{\ell \in V_n} \frac{\E \bigl| \eta_{\ell,n} - \E \{ \eta_{\ell,n} \} \bigr|^p}{\Bigl( \text{Var}\bigl\{ T_n^{(K)} \bigr\} \Bigr)^{p/2}} \to 0 \ \ \text{as } n \to \infty\,.
\end{equation}
Since the proof for showing this varies depending on the limit of
$nf(R_ne_1)$, we divide the argument into three different cases.
Suppose first that $nf(R_ne_1) \to 0$ as $n \to \infty$. Let
$Z_{\ell,n}$ denote the number of points in $\Pn$ lying in
$$
\text{Tube} (Q_\ell; kt_m) := \bigl\{ x \in \bbr^d: \inf_{y \in Q_{\ell}} ||x-y|| \leq kt_m \bigr\}\,.
$$
Then, $Z_{\ell,n}$ has a Poisson distribution with mean
$n\int_{\text{Tube}(Q_\ell; kt_m)}f(z)dz$. Using Potter's bound,
we see that $Z_{\ell,n}$ is stochastically dominated by another
Poisson random variable $Z_n$ with mean $C^* nf(R_ne_1)$.
Observing that
$$
|\eta_{\ell,n}| \leq C^* \begin{pmatrix} Z_{\ell,n}  \\ k \end{pmatrix}\,,
$$
we have, for $q=1,2,3,4$,
$$
\E|\eta_{\ell,n}|^q \leq C^* \E \begin{pmatrix} Z_{\ell,n}  \\ k \end{pmatrix}^q \leq C^* \E \begin{pmatrix} Z_n  \\ k \end{pmatrix}^q \leq C^* \bigl( nf(R_ne_1) \bigr)^k,
$$
where in the last step we used the assumption $nf(R_ne_1) \to 0$. \\
It now follows that for $p=3,4$,
$$
\max_{\ell \in V_n} E \bigl| \eta_{\ell,n} -\E \{ \eta_{\ell,n} \} \bigr|^p \leq C^* \bigl( nf(R_ne_1) \bigr)^k.
$$
Therefore,
\begin{align*}
R_n^d\, \max_{\ell \in V_n} \frac{\E \bigl| \eta_{\ell,n} - \E \{ \eta_{\ell,n} \} \bigr|^p}{\Bigl( \text{Var}\bigl\{ T_n^{(K)} \bigr\} \Bigr)^{p/2}} &\leq C^* R_n^d \frac{\bigl( nf(R_ne_1) \bigr)^k}{ \bigl( n^k R_n^d f(R_ne_1)^k \gamma_K \bigr)^{p/2}}  \\
&= \frac{C^*}{\gamma_K^{p/2}} \bigl( n^k R_n^d f(R_ne_1)^k \bigr)^{1-p/2} \to0\,, \ \ \ n \to\infty\,,
\end{align*}
where the last convergence follows from \eqref{e:normalizing.heavy}.

In the case of $nf(R_ne_1) \to \xi \in (0,\infty)$, the argument
for proving \eqref{e:normal.approx} is very similar to, or even
easier than, the previous case, so we omit it.

Finally, suppose that $nf(R_ne_1) \to \infty$ as $n \to \infty$.
We begin by establishing an appropriate upper bound for the fourth
moment expectation
\begin{equation}  \label{e:bino.expansion}
\E \bigl| \eta_{\ell,n} - \E \{ \eta_{\ell,n} \} \bigr|^4 = \sum_{j=0}^4 \begin{pmatrix} 4 \\ j \end{pmatrix} (-1)^j \E \{ \eta_{\ell,n}^j \} \bigl( \E \{ \eta_{\ell,n} \} \bigr)^{4-j}.
\end{equation}
Letting
$$
g_{\ell,n}(\Y) := \sum_{j=1}^m a_j h_{t_j}(\Y)\, \one \bigl\{ m(\Y) \geq R_n, \ \text{Max}(\Y) \in  \text{Ann} (R_n,KR_n) \cap Q_{\ell} \bigr\}\,,
$$
we see that for every $j \in \{0,\dots,4\}$,
$$
F_n(j):=\E \{ \eta_{\ell,n}^j \}\bigl( \E \{ \eta_{\ell,n} \} \bigr)^{4-j}
$$
can be denoted as the expectation of a quadruple sum
\begin{equation}  \label{e:quadruple}
\E \left\{ \sum_{\Y_1 \subset \Pn^{(1)}} \sum_{\Y_2 \subset \Pn^{(2)}} \sum_{\Y_3 \subset \Pn^{(3)}} \sum_{\Y_4 \subset \Pn^{(4)}} g_{\ell,n}(\Y_1)\, g_{\ell,n}(\Y_2)\, g_{\ell,n}(\Y_3)\, g_{\ell,n}(\Y_4)\, \right\},
\end{equation}
where each of $\Pn^{(1)}, \dots, \Pn^{(4)}$ is either equal to or
an independent copy of one of the others. By definition, each
$\Y_i$ is a finite collection of $d$-dimensional vectors. If, in
particular, $|\Y_1 \cup \Y_2 \cup \Y_3 \cup \Y_4| = 4k$, i.e., any
two of $\Y_i$, $i=1,\dots,4$ have no common elements, then the
Palm theory given in the Appendix reveals that \eqref{e:quadruple}
is equal to $\bigl( \E \{ \eta_{\ell,n} \} \bigr)^4$. Then, in
this case, their overall contribution to \eqref{e:bino.expansion} is
identically zero, because
$$
\sum_{j=0}^4  \begin{pmatrix} 4 \\ j \end{pmatrix} (-1)^j \bigl( \E \{ \eta_{\ell,n} \} \bigr)^{4} = 0\,.
$$

Next, suppose that $|\Y_1 \cup \Y_2 \cup \Y_3 \cup \Y_4| = 4k-1$,
i.e., there is a pair $(\Y_i, \Y_j)$, $i \neq j$ having exactly
one element in common and no other common elements between
$\Y_i$'s are present. In this case, \eqref{e:quadruple} can be
written as
\begin{equation}  \label{e:second.order.term}
\frac{n^{2k-1}}{\bigl( (k-1)! \bigr)^2}\, \E \Bigl\{ g_{\ell,n}(\Y_1)\, g_{\ell,n}(\Y_2)\, \one \bigl\{  |\Y_1 \cap \Y_2|=1 \bigr\} \Bigr\}
 \left( \frac{n^k}{k!} \E \bigl\{ g_{\ell,n}(\Y) \bigr\} \right)^2.
\end{equation}
In particular, \eqref{e:second.order.term} appears once in
$F_n(2)$, $\tiny  \begin{pmatrix} 3 \\ 2 \end{pmatrix} $ times in
$F_n(3)$, and $\tiny  \begin{pmatrix} 4 \\ 2 \end{pmatrix} $ times
in $F_n(4)$. Thus, the total contribution to
\eqref{e:bino.expansion} sums up to
$$
\biggl\{ \begin{pmatrix}  4 \\ 2 \end{pmatrix} (-1)^2 + \begin{pmatrix}  4 \\ 3 \end{pmatrix} (-1)^3  \begin{pmatrix}  3 \\ 2 \end{pmatrix}  + \begin{pmatrix}  4 \\ 4 \end{pmatrix} (-1)^4  \begin{pmatrix}  4 \\ 2 \end{pmatrix}  \biggr\} \times \eqref{e:second.order.term} = 0\,.
$$

We may assume, therefore, that $|\Y_1 \cup \Y_2 \cup \Y_3 \cup
\Y_4| \leq 4k-2$. Let us start with $|\Y_1 \cup \Y_2 \cup \Y_3
\cup \Y_4| = 4k-2$, where we shall  examine in particular the case
in which $\Pn^{(1)}=\Pn^{(2)}=\Pn^{(3)}=\Pn^{(4)}$, $|\Y_1 \cap
\Y_2|=2$ and no other common elements  between $\Y_i$'s exist. The
argument for the other cases will be omitted because they can be
handled in the same manner. Then, by Palm theory,
\eqref{e:quadruple} is equal to
\begin{equation}  \label{e:2overlaps}
\frac{n^{2k-2}}{2 \bigl( (k-2)! \bigr)^2}\, \E \Bigl\{ g_{\ell,n}(\Y_1)\, g_{\ell,n}(\Y_2)\, \one \bigl\{ |\Y_1 \cap \Y_2| = 2 \bigr\} \Bigr\} \left( \frac{n^k}{k!} \E \bigl\{ g_{\ell,n}(\Y) \bigr\} \right)^2.
\end{equation}
Because of Potter's bound, together with the fact that $Q_\ell$
intersects with $\text{Ann}(R_n,KR_n)$,
\begin{align*}
\Bigl|  \E &\Bigl\{ g_{\ell,n}(Y_1)\, g_{\ell,n}(\Y_2)\, \one \bigl\{ |\Y_1 \cap \Y_2| = 2 \bigr\} \Bigr\} \Bigr|  \\
&\leq C^* \Bigl( \P \bigl\{ X_1 \in \text{Tube} (Q_\ell; kt_m) \bigr\} \Bigr)^{2k-2} \leq C^* f(R_ne_1)^{2k-2}.
\end{align*}
Similarly, we can obtain
$$
\Bigl| \E \bigl\{ g_{\ell,n}(\Y) \bigr\} \Bigr| \leq C^* f(R_ne_1)^k,
$$
and therefore, the absolute value of \eqref{e:2overlaps},
equivalently that of \eqref{e:quadruple}, is bounded above by $C^*
\bigl( nf(R_ne_1) \bigr)^{4k-2} $.

A similar argument proves that if $|\Y_1 \cup \Y_2 \cup \Y_3 \cup
\Y_4| = 4k-q$ for some $q \geq 3$, the absolute value of
\eqref{e:quadruple} is bounded above by $C^*\bigl( nf(R_ne_1)
\bigr)^{4k-q}$. Putting these facts altogether, while recalling
$nf(R_ne_1) \to \infty$ as $n\to \infty$, we may conclude that
$$
\E \bigl| \eta_{\ell,n} - \E \{ \eta_{\ell,n} \} \bigr|^4 \leq C^* \bigl( nf(R_ne_1) \bigr)^{4k-2}.
$$
Now, it is easy to check \eqref{e:normal.approx}.

In terms of the third moment expectation $\E \bigl| \eta_{\ell,n}
- \E \{ \eta_{\ell,n} \} \bigr|^3$, we apply H\"{o}lder's
inequality to obtain
$$
\E \bigl| \eta_{\ell,n} - \E \{ \eta_{\ell,n} \} \bigr|^3 \leq \Bigl( \E \bigl| \eta_{\ell,n} - \E \{ \eta_{\ell,n} \} \bigr|^4  \Bigr)^{3/4} \leq C^*  \bigl( nf(R_ne_1) \bigr)^{3k-3/2}.
$$
Again, it is easy to prove \eqref{e:normal.approx}.

Now, we have obtained a CLT in \eqref{e:int.CLT} as required, regardless of the limit of $nf(R_ne_1)$.
\end{proof}
An important claim is that once the tightness of each $\BX_n^+$
and $\BX_n^-$ is established in the space $\mathcal D [0,\infty)$ which is equipped with the Skorohod $J_1$-topology, the proof of Theorem
\ref{t:main.heavy} is complete. To see this, suppose that
$\BX_n^+$ and $\BX_n^-$ were both tight in $\mathcal D
[0,\infty)$. Then, a joint process $(\BX_n^+, \BX_n^-)$ is tight
as well in $\mathcal D [0,\infty) \times \mathcal D [0,\infty)$,
which is endowed with the product topology. Because of the already
established finite-dimensional weak convergence of $(\BX_n^+,
\BX_n^-)$, every subsequential limit of $(\BX_n^+, \BX_n^-)$
coincides with the limiting process in Proposition
\ref{p:fidi.heavy}. This in turn implies the weak convergence of
$(\BX_n^+, \BX_n^-)$ in $\mathcal D [0,\infty) \times \mathcal D
[0,\infty)$. Using the basic fact that the map $(x,y) \to x-y$
from $\mathcal D [0,\infty) \times \mathcal D [0,\infty)$ to
$\mathcal D [0,\infty) $ is continuous at $(x,y) \in \mathcal
C[0,\infty) \times \mathcal C[0,\infty)$, while recalling that the
limits in Proposition \ref{p:fidi.heavy} all have continuous
sample paths, the continuous mapping theorem gives weak
convergence of $\BX_n = \BX_n^+ - \BX_n^-$ in $\mathcal D
[0,\infty)$.
\begin{proposition}  \label{p:tightness.heavy}
The sequences $(\BX_n^+)$ and $(\BX_n^-)$ are both tight in
$\mathcal D [0,\infty)$, irrespective of the limit of
$nf(R_ne_1)$.
\end{proposition}
\begin{proof}
We prove the tightness of $(\BX_n^+)$ only, in the space $\mathcal D [0,L]$ for any fixed $L>0$. For notational ease, however, we omit the superscript ``+" from all the functions and objects during the proof. By Theorem 13.5 of \cite{billingsley:1999}, it is
sufficient to show that there exists $B>0$ such that
$$
\E \Bigl\{  \bigl( X_n(t) - X_n(s) \bigr)^2 \bigl( X_n(s) - X_n(r) \bigr)^2 \Bigr\} \leq B (t-r)^2
$$
for all $0 \leq r \leq s \leq t \leq L$ and $n \geq 1$.

For typographical convenience, we use shorthand notations
\eqref{e:def.hts}, \eqref{e:def.hnts}, and further,
\begin{align*}
\xi_{n,t,s} &:= \sum_{\Y \subset \Pn} h_{n, t,s}(\Y)\,.
\end{align*}
Then,
\begin{align*}
\E \Bigl\{  \bigl( X_n(t) - X_n(s) \bigr)^2 \bigl( X_n(s) - X_n(r) \bigr)^2 \Bigr\}
&= \tau_n^{-2} \E \Bigl\{  \bigl( \xi_{n,t,s} - \E \{ \xi_{n,t,s} \} \bigr)^2 \bigl( \xi_{n,s,r} - \E \{ \xi_{n,s,r} \} \bigr)^2 \Bigr\} \\
&= \tau_n^{-2} \sum_{p=0}^2 \sum_{q=0}^2 \begin{pmatrix} 2 \\ p \end{pmatrix} \begin{pmatrix} 2 \\ q \end{pmatrix} (-1)^{p+q} F_n(p,q)\,,
\end{align*}
where
$$
F_n(p,q) = \E \{ \xi_{n,t,s}^p \xi_{n,s,r}^q \} \bigl( \E \{ \xi_{n,t,s} \} \bigr)^{2-p} \bigl( \E \{ \xi_{n,s,r} \} \bigr)^{2-q}.
$$
Note that for every $p, q \in \{ 0,1,2 \}$, $F_n(p,q)$ can be represented by
\begin{equation}  \label{e:quadruple1}
\E \left\{ \sum_{\Y_1 \subset \Pn^{(1)}} \sum_{\Y_2 \subset \Pn^{(2)}} \sum_{\Y_3 \subset \Pn^{(3)}} \sum_{\Y_4 \subset \Pn^{(4)}} h_{n,t,s}(\Y_1)\, h_{n,t,s}(\Y_2)\, h_{n,s,r}(\Y_3)\, h_{n,s,r}(\Y_4)\, \right\},
\end{equation}
where each of $\Pn^{(1)}, \dots, \Pn^{(4)}$ is either equal to or an independent copy of one of the others.

According to the Palm theory given in the Appendix, if $|\Y_1 \cup
\Y_2 \cup \Y_3 \cup \Y_4 | = 4k$, i.e., any two of $\Y_i$ have no
common elements, then \eqref{e:quadruple1} reduces to $\bigl( \E
\{ \xi_{n,t,s} \} \bigr)^2 \bigl( \E \{ \xi_{n,s,r} \} \bigr)^2$.
Then, an overall contribution in this case identically vanishes,
since
$$
\sum_{p=0}^2 \sum_{q=0}^2 \begin{pmatrix} 2 \\ p \end{pmatrix} \begin{pmatrix} 2 \\ q \end{pmatrix} (-1)^{p+q} \bigl( \E \{ \xi_{n,t,s} \} \bigr)^2 \bigl( \E \{ \xi_{n,s,r} \} \bigr)^2 = 0\,.
$$

In the following, we examine the case in which at least one common
element exists between $\Y_i$'s. First, for $\ell = 1,\dots,k$, we
count the number of times
\begin{equation}  \label{e:pattern1}
\E \Bigl\{ \sum_{\Y_1 \subset \Pn}\sum_{\Y_2 \subset \Pn}h_{n,t,s}(\Y_1)\, h_{n,t,s}(\Y_2)\, \one \bigl\{  |\Y_1 \cap \Y_2|=\ell \bigr\} \Bigr\}
 \bigl(\E \{ \xi_{n,s,r} \} \bigr)^2
\end{equation}
appears in each $F_n(p,q)$. Indeed, \eqref{e:pattern1} appears
only once  in $F_n(2,0)$, $F_n(2,1)$, and $F_n(2,2)$. Therefore,
the total contribution amounts to
$$
\left[ \begin{pmatrix} 2 \\ 2 \end{pmatrix}\begin{pmatrix} 2 \\ 0 \end{pmatrix}(-1)^{2+0} + \begin{pmatrix} 2 \\ 2 \end{pmatrix}\begin{pmatrix} 2 \\ 1 \end{pmatrix}(-1)^{2+1} + \begin{pmatrix} 2 \\ 2 \end{pmatrix}\begin{pmatrix} 2 \\ 2 \end{pmatrix}(-1)^{2+2} \right] \times \eqref{e:pattern1} = 0\,.
$$
Similarly, for every $\ell = 1,\dots, k$, no contribution is made
by
$$
\E \Bigl\{ \sum_{\Y_1 \subset \Pn}\sum_{\Y_2 \subset \Pn}h_{n,s,r}(\Y_1)\, h_{n,s,r}(\Y_2)\, \one \bigl\{  |\Y_1 \cap \Y_2|=\ell \bigr\} \Bigr\}
 \bigl(\E \{ \xi_{n,t,s} \} \bigr)^2.
$$

Subsequently, for $\ell=1,\dots,k$, we explore the presence of
\begin{equation}  \label{e:pattern2}
\E \Bigl\{ \sum_{\Y_1 \subset \Pn}\sum_{\Y_2 \subset \Pn}h_{n,t,s}(\Y_1)\, h_{n,s,r}(\Y_2)\, \one \bigl\{  |\Y_1 \cap \Y_2|=\ell \bigr\} \Bigr\}
 \E \{ \xi_{n,t,s} \} \E \{ \xi_{n,s,r} \}\,.
\end{equation}
One can immediately check that \eqref{e:pattern2} appears once in
$F_n(1,1)$, twice in $F_n(2,1)$, twice in $F_n(1,2)$, and four
times in $F_n(2,2)$. However, their total contribution disappears
again, because
\begin{multline*}
\biggl[ \begin{pmatrix} 2 \\ 1 \end{pmatrix}\begin{pmatrix} 2 \\ 1 \end{pmatrix}(-1)^{1+1} + \begin{pmatrix} 2 \\ 2 \end{pmatrix}\begin{pmatrix} 2 \\ 1 \end{pmatrix}(-1)^{2+1} \cdot 2 \\ + \begin{pmatrix} 2 \\ 1 \end{pmatrix}\begin{pmatrix} 2 \\ 2 \end{pmatrix}(-1)^{1+2} \cdot 2
+ \begin{pmatrix} 2 \\ 2 \end{pmatrix}\begin{pmatrix} 2 \\ 2 \end{pmatrix}(-1)^{2+2} \cdot 4 \biggr] \times \eqref{e:pattern2} = 0\,.
\end{multline*}

Next, let $\ell_i \in \{ 0,\dots,k \}$, $i=1,2,3$, $\ell \in \{
2,\dots,2k\}$ such that at least two of $\ell_i$'s are non-zero,
so that we should examine the appearance of
\begin{multline}
\E \Bigl\{  \sum_{\Y_1 \subset \Pn}\sum_{\Y_2 \subset \Pn} \sum_{\Y_3 \subset \Pn} h_{n,t,s}(\Y_1)\, h_{n,t,s}(\Y_2)\,h_{n,s,r}(\Y_3)\,  \label{e:pattern3} \\
\times \one \bigl\{ |\Y_1 \cap \Y_2 | =\ell_1, \,  |\Y_1 \cap \Y_3 | =\ell_2, \, |\Y_2 \cap \Y_3 | =\ell_3, \, |\Y_1 \cup \Y_2 \cup \Y_3| = 3k-\ell \bigr\} \Bigr\}\, \E \{ \xi_{n,s,r} \}\,.
\end{multline}
This actually appears once in $F_n(2,1)$ and twice in $F_n(2,2)$; therefore, their overall contribution is
$$
\biggl[ \begin{pmatrix} 2 \\ 2 \end{pmatrix}\begin{pmatrix} 2 \\ 1 \end{pmatrix}(-1)^{2+1} + \begin{pmatrix} 2 \\ 2 \end{pmatrix}\begin{pmatrix} 2 \\ 2 \end{pmatrix}(-1)^{2+2} \cdot 2  \biggr] \times \eqref{e:pattern3} = 0\,.
$$
For the same reason, we can ignore the presence of
\begin{multline*}
\E \Bigl\{  \sum_{\Y_1 \subset \Pn}\sum_{\Y_2 \subset \Pn} \sum_{\Y_3 \subset \Pn} h_{n,t,s}(\Y_1)\, h_{n,s,r}(\Y_2)\,h_{n,s,r}(\Y_3)\,  \\
\times \one \bigl\{ |\Y_1 \cap \Y_2 | =\ell_1, \,  |\Y_1 \cap \Y_3 | =\ell_2, \, |\Y_2 \cap \Y_3 | =\ell_3, \, |\Y_1 \cup \Y_2 \cup \Y_3| = 3k-\ell \bigr\} \Bigr\}\, \E \{ \xi_{n,t,s} \}\,.
\end{multline*}
where $\ell_i \in \{ 0,\dots,k \}$, $i=1,2,3$, $\ell \in \{
2,\dots,2k\}$ such that at least two of $\ell_i$'s are non-zero.

Putting these calculations altogether, we find that the tightness
follows, once we can show that there exists $B>0$ such that
\begin{multline}
\tau_n^{-2} \E \Bigl\{ \sum_{\Y_1 \subset \Pn}\sum_{\Y_2 \subset \Pn} \sum_{\Y_3 \subset \Pn}\sum_{\Y_4 \subset \Pn} h_{n,t,s}(\Y_1)\, h_{n,t,s}(\Y_2)\,h_{n,s,r}(\Y_3)\,h_{n,s,r}(\Y_4)\, \label{e:pattern4} \\
\times \one \bigl\{ \text{each } \Y_i \text{ has at least one common elements with}  \\
\text{at least one of the other three} \bigr\} \Bigr\} \leq B(t-r)^2
\end{multline}
for all $0 \leq r \leq s \leq t \leq L$ and $n \geq 1$. We  need to check only the following possibilities. \\
$[\text{I}]$ $\ell := |\Y_1 \cap \Y_2| \in \{ 1,\dots,k \}$, $\ellp := |\Y_3 \cap \Y_4| \in \{ 1,\dots,k \}$,
and $(\Y_1 \cup \Y_2) \cap (\Y_3 \cup \Y_4) = \emptyset$. \\
$[\text{I} \hspace{-1pt }\text{I}]$ $\ell := |\Y_2 \cap \Y_3| \in \{ 1,\dots,k \}$, $\ellp := |\Y_1 \cap \Y_4| \in \{ 1,\dots,k \}$, and $(\Y_2 \cup \Y_3) \cap (\Y_1 \cup \Y_4) = \emptyset$. \\
$[\text{I}\hspace{-1pt}\text{I}\hspace{-1pt}\text{I}]$. Each
$\Y_i$ has at least one common element with at least one of the
other three, but neither $[\text{I}]$ or $[\text{I}\hspace{-1pt
}\text{I}]$ is true.

For example, if $|\Y_1 \cap \Y_2|=2$, $|\Y_1 \cap \Y_3|=3$, $|\Y_2
\cap \Y_4|=1$, and there are no other common elements between
$\Y_i$'s, then it falls into category
$[\text{I}\hspace{-1pt}\text{I}\hspace{-1pt}\text{I}]$, where,
unlike $[\text{I}]$ or $[\text{I} \hspace{-1pt }\text{I}]$, the
expectation in \eqref{e:pattern4} can no longer be separated by
the Palm theory.

Denoting by $A$ the left-hand side of \eqref{e:pattern4},  let us
start with case $[\text{I}]$. As a result of Palm theory,
\begin{align*}
A &= \tau_n^{-1} \frac{n^{2k-\ell}}{\ell ! \bigl( (k-\ell)! \bigr)^2}\, \E \Bigl\{ h_{n,t,s}(\Y_1)\, h_{n,t,s}(\Y_2)\, \one \bigl\{ |\Y_1 \cap \Y_2| = \ell \bigr\} \Bigr\} \\
&\quad \times \tau_n^{-1} \frac{n^{2k-\ellp}}{\ellp ! \bigl( (k-\ellp)! \bigr)^2}\, \E \Bigl\{ h_{n,s,r}(\Y_3)\, h_{n,s,r}(\Y_4)\, \one \bigl\{ |\Y_3 \cap \Y_4| = \ellp \bigr\} \Bigr\} \\
&:= A_1 \times A_2\,.
\end{align*}
Proceeding as in the calculation of Proposition
\ref{p:cova.heavy}, we obtain
\begin{align}
A_1 &\leq C^* \tau_n^{-1} n^{2k-\ell} R_n^d f(R_ne_1)^{2k-\ell} \int_{(\bbr^d)^{\ell-1}} \hspace{-10pt} d\by \int_{(\bbr^d)^{k-\ell}} \hspace{-10pt} d\bz_2 \int_{(\bbr^d)^{k-\ell}} \hspace{-10pt} d\bz_1\, h_{t,s}(0,\by, \bz_1)\, h_{t,s}(0,\by, \bz_2)\,,  \label{e:upper.E1}\\
A_2 &\leq C^* \tau_n^{-1} n^{2k-\ellp} R_n^d f(R_ne_1)^{2k-\ellp} \int_{(\bbr^d)^{\ellp-1}} \hspace{-10pt} d\by \int_{(\bbr^d)^{k-\ellp}} \hspace{-10pt} d\bz_2 \int_{(\bbr^d)^{k-\ellp}} \hspace{-10pt} d\bz_1\, h_{s,r}(0,\by, \bz_1)\, h_{s,r}(0,\by, \bz_2)\,.  \label{e:upper.E2}
\end{align}
Notice that $h_t$ is increasing in $t$ in the sense of \eqref{e:ind.increase+} (recall that the superscript ``+" is suppressed during the proof). 
It also follows from \eqref{e:close.enough.decomp.dyna} that the triple integral in
\eqref{e:upper.E1} is unchanged if the integral domain is
restricted to $\bigl( B(0,kL) \bigr)^{\ell-1} \times \bigl(
B(0,kL) \bigr)^{k-\ell} \times \bigl( B(0,kL) \bigr)^{k-\ell}$.
 Therefore, with $\lambda$ being the
Lebesgue measure on $(\bbr^d)^{k-\ell}$,
\begin{align*}
\int_{(\bbr^d)^{\ell-1}} \hspace{-10pt} d\by &\int_{(\bbr^d)^{k-\ell}} \hspace{-10pt} d\bz_2 \int_{(\bbr^d)^{k-\ell}} \hspace{-10pt} d\bz_1\, h_{t,s}(0,\by, \bz_1)\, h_{t,s}(0,\by, \bz_2) \\
&\leq \lambda \bigl\{ \bigl( B(0,kL) \bigr)^{k-\ell} \bigr\}\, \int_{(\bbr^d)^{\ell-1}}\int_{(\bbr^d)^{k-\ell}} h_{t,s}(0,\by, \bz) d\by d\bz \\
&= \lambda \bigl\{ \bigl( B(0,kL) \bigr)^{k-\ell} \bigr\}\, \bigl( t^{d(k-1)} - s^{d(k-1)} \bigr) \int_{(\bbr^d)^{k-1}} h(0,\by) d\by \\
&\leq C^*(t-r)\,.
\end{align*}
Applying the same manipulation to the triple integral in
\eqref{e:upper.E2}, we obtain
$$
A \leq C^* \tau_n^{-2} n^{4k-\ell-\ellp} R_n^{2d} f(R_ne_1)^{4k-\ell-\ellp} (t-r)^2.
$$
It remains to check that $\sup_n \tau_n^{-2} n^{4k-\ell-\ellp}
R_n^{2d} f(R_ne_1)^{4k-\ell-\ellp}  < \infty$, which is, however,
easy to prove, irrespective of the definition of $\tau_n$. Now
case $[\text{I}]$ is done.

Next, we turn to case $[\text{I}\hspace{-1pt}\text{I}]$. As a
consequence of the same operation as in $[\text{I}]$, we obtain
the same upper bound for $A$ up to multiplicative constants.

Finally, we proceed to case
$[\text{I}\hspace{-1pt}\text{I}\hspace{-1pt}\text{I}]$. Let $\ell
:= 4k-|\Y_1 \cup \Y_2 \cup \Y_3 \cup \Y_4|$; then, it must be that
$3 \leq \ell \leq 3k$. It follows from Palm theory that
$$
A = C^* \tau_n^{-2} n^{4k-\ell} \E \bigl\{ h_{n,t,s}(\Y_1)\, h_{n,t,s}(\Y_2)\,h_{n,s,r}(\Y_3)\,h_{n,s,r}(\Y_4) \bigr\}
$$
with $(\Y_1,\dots,\Y_4)$ satisfying requirements in case
$[\text{I}\hspace{-1pt}\text{I}\hspace{-1pt}\text{I}]$. In
particular, $(\Y_1 \cup \Y_2) \cap (\Y_3 \cup \Y_4)$ must be
non-empty; hence, we may assume without loss of generality that
$\Y_1 \cap \Y_3\neq \emptyset$. Set $\ellp := |\Y_1 \cap \Y_3| \in
\{ 1,\dots,k \}$. By \eqref{e:ind.increase+} and \eqref{e:close.enough.decomp.dyna}, we have
$$
A \leq C^* \tau_n^{-2} n^{4k-\ell} R_n^d f(R_ne_1)^{4k-\ell} \int_{(\bbr^d)^{\ellp-1}} \hspace{-10pt} d\by \int_{(\bbr^d)^{k-\ellp}} \hspace{-10pt} d\bz_2 \int_{(\bbr^d)^{k-\ellp}} d\bz_1\, h_{t,s}(0,\by, \bz_1)\, h_{t,s}(0,\by, \bz_2)\,,
$$
Because of Lemma \ref{l:used.for.tightness},
\begin{align*}
A \leq C^* \tau_n^{-2} n^{4k-\ell} R_n^d f(R_ne_1)^{4k-\ell} (t-r)^2.
\end{align*}
Once again, verifying
$$
\sup_n \tau_n^{-2} n^{4k-\ell} R_n^d f(R_ne_1)^{4k-\ell} < \infty
$$
is elementary, and hence, we have completed the proof of
\eqref{e:pattern4} as required.
\end{proof}

\subsection{Exponentially Decaying Tail Case}  \label{s:proof.light}

We start by defining a subgraph counting process with restricted
domain. For $0 \leq K < L \leq \infty$, we define
\begin{align*}
G_{n,K,L}(t) &= \sum_{\Y \subset \mathcal{P}_n} h_t (\Y)\, \one \bigl\{ m(\Y) \geq R_n\,, \ a(R_n)^{-1}\bigl(\text{Max}(\Y) - R_n \bigr) \in [K,L) \bigr\} \\
&:= \sum_{\Y \subset \mathcal{P}_n} h_{n,t,K,L} (\Y)\,,
\end{align*}
and
\begin{align*}
G_{n,K,L}^{\pm}(t) &= \sum_{\Y \subset \mathcal{P}_n} h_t^{\pm} (\Y)\, \one \bigl\{ m(\Y) \geq R_n\,, \  a(R_n)^{-1}\bigl(\text{Max}(\Y) - R_n \bigr) \in [K,L) \bigr\} \\
&:= \sum_{\Y \subset \mathcal{P}_n} h_{n,t,K,L}^{\pm} (\Y)\,,
\end{align*}
where  $(R_n)$ satisfies \eqref{e:normalizing.light}. For the
special case $K=0$ and $L=\infty$, we denote $G_n(t) =
G_{n,0,\infty}(t)$ and $G_n^{\pm}(t) = G_{n,0,\infty}^{\pm}(t)$.
The centered and scaled versions of the subgraph counting process
are
\begin{align}
X_n(t) &= \tau_n^{-1/2} \Bigl( G_n(t) - \E\bigl\{ G_n(t) \bigr\} \Bigr)\,, \label{e:def.Xn}\\
X_n^{\pm}(t) &= \tau_n^{-1/2} \Bigl( G_n^{\pm}(t) - \E\bigl\{ G_n^{\pm}(t) \bigr\} \Bigr) \label{e:def.Xnpm}\,,
\end{align}
where $(\tau_n)$ is given in \eqref{e:tau.light}. As seen in the
regularly varying tail case, we first need to know the growing
rate of the covariances of $G_{n,K,L}(t)$. Before presenting the
results, we introduce for $\ell = 1,\dots,k$,
\begin{align*}
M_{\ell, K, L}(t,s) := D_\ell &\int_0^{\infty} \int_{(\bbr^d)^{2k-\ell-1}} \hspace{-10pt}e^{ -(2k-\ell)\rho - c^{-1} \sum_{i=1}^{2k-\ell-1} \langle e_1,y_i \rangle }\, \\
&\times  \one \bigl\{  \by \in E_{K,L}^{(\ell)} (\rho,e_1) \bigr\}\, h_{t,s}^{(\ell)}(0,\by)\, d\by d\rho\,, \ \ t,s\geq 0\,,
\end{align*}
where $D_\ell $ is given in \eqref{e:def.D.ell}, $h_{t,s}^{(\ell)}(0,\by)$ is defined in \eqref{e:def.h.ell}, and for $\rho > 0$ and $\theta \in S_{d-1}$,
\begin{align*}
E_{K,L}^{(\ell)}(\rho, \theta) = \Bigl\{  \by \in (\bbr^d)^{2k-\ell-1}:\,  &\rho + c^{-1} \langle \theta, y_i \rangle \geq 0\,, \ i=1,\dots, 2k-\ell-1\,, \\
&K \leq \max \bigl\{ \rho\,, \rho + c^{-1} \hspace{-5pt} \max_{i = 1,\dots, k-1}\langle \theta, y_i \rangle  \bigr\} < L\,, \\
&K \leq \max \bigl\{ \rho\,, \rho + c^{-1} \hspace{-5pt} \max_{i=1,\dots,\ell-1,k,\dots,2k-\ell-1}\langle \theta, y_i \rangle  \bigr\} < L \, \Bigr\}\,.
\end{align*}
Note that $M_{\ell, 0, \infty}(t,s)$ completely matches \eqref{e:cov.comp.light}.

\begin{proposition}  \label{p:cova.light}
Assume the conditions of Theorem \ref{t:main.light}. Let $0 \leq K < L \leq \infty$. \\
$(i)$ If $nf(R_ne_1) \to 0$ as $n \to \infty$, then
$$
\tau_n^{-1} \text{Cov} \bigl( G_{n,K,L}(t), G_{n,K,L}(s) \bigr) \to M_{k,K,L}(t,s)\,, \ \ \ n\to \infty\,.
$$
$(ii)$ If $nf(R_ne_1) \to \xi \in (0,\infty)$ as $n \to \infty$, then
$$
\tau_n^{-1} \text{Cov} \bigl( G_{n,K,L}(t), G_{n,K,L}(s) \bigr) \to \sum_{\ell=1}^k  \xi^{2k-\ell}M_{\ell, K, L}(t,s)\,, \ \ \ n\to \infty\,.
$$
$(iii)$ If $nf(R_ne_1) \to \infty$ as $n \to \infty$, then
$$
\tau_n^{-1} \text{Cov} \bigl( G_{n,K,L}(t), G_{n,K,L}(s) \bigr) \to M_{1,K,L}(t,s)\,, \ \ \ n\to \infty\,.
$$
\end{proposition}
\begin{proof}
As argued in Proposition \ref{p:cova.heavy}, with the multiple
applications of Palm theory, one can write
$$
\text{Cov} \bigl( G_{n,K,L}(t)\,, G_{n,K,L}(s) \bigr)
= \sum_{\ell=1}^k \frac{n^{2k-\ell}}{\ell ! \bigl( (k-\ell)! \bigr)^2}\, \E \Bigl\{h_{n,t, K, L} (\Y_1)\, h_{n,s, K, L} (\Y_2)\, \one \bigl\{ \, |\Y_1 \cap \Y_2| = \ell\, \bigr\} \Bigr\}.
$$
Define for $\ell \in \{ 1,\dots,k \}$,
\begin{multline*}
F_n^{(\ell)} (K,L) := \bigl\{ \bx \in (\bbr^d)^{2k-\ell}: a(R_n)^{-1} \bigl( \text{Max}(x_1,\dots, x_k) - R_n \bigr) \in [K,L)\,, \\
 a(R_n)^{-1} \bigl( \text{Max}(x_1,\dots, x_\ell, x_{k+1}, \dots, x_{2k-\ell}) - R_n \bigr) \in [K,L) \bigr\}\,.
\end{multline*}
By the change of variables $\bx \rightarrow (x,x+\by)$ with $\bx
\in (\bbr^d)^{2k-\ell}$, $x \in \bbr^d$, $\by \in
(\bbr^d)^{2k-\ell-1}$, together with invariance
\eqref{e:location.inv},
\begin{align*}
\E \Bigl\{ &h_{n,t,K,L}(\Y_1)\, h_{n,s,K,L}(\Y_2)\, \one \bigl\{ \, |\Y_1 \cap \Y_2| = \ell\, \bigr\} \Bigr\} \\
&= \int_{(\bbr^d)^{2k-\ell}} f(\bx)\, \one \bigl\{ m(\bx) \geq R_n \bigr\}\,  h_{t,s}^{(\ell)} (\bx)\,\one \bigl\{ \bx \in F_n^{(\ell)}(K,L) \bigr\} d\bx  \\
&= \int_{\bbr^d} \int_{(\bbr^d)^{2k-\ell-1}} f(x)\, f(x+\by)\, \one \bigl\{ m(x, x + \by) \geq R_n \bigr\}\, h_{t,s}^{(\ell)} (0, \by)\,  \\
&\quad \times \one \bigl\{ (x, x+\by) \in F_n^{(\ell)}(K,L) \bigr\} d\by dx\,.
\end{align*}
Let $J_k$ denote the last integral. Further calculation by the
polar coordinate transform $x \to (r,\theta)$ with $J(\theta) =
|\partial x / \partial \theta|$ and the change of variable $\rho =
a(R_n)^{-1} (r-R_n) $ yields
\begin{align}
J_k = &a(R_{n}) R_{n}^{d-1} f(R_{n}e_1)^{2k-\ell} \int_{S^{d-1}} \hspace{-7pt} J(\theta)d\theta \int_0^{\infty} d\rho \int_{(\bbr^d)^{2k-\ell-1}} d\by\,  \label{e:Jk}\\
&\times \left( 1 + \frac{a(R_{n})}{R_{n}} \rho \right)^{d-1} \frac{f\Bigl( \bigl( R_{n} + a(R_{n}) \rho \bigr) e_1 \Bigr)}{f(R_ne_1)} \notag \\
&\times \prod_{i=1}^{2k-\ell-1} f(R_{n}e_1)^{-1} f\Bigl( \|(R_{n}+a(R_{n})\rho)\theta + y_i\|e_1 \Bigr)\, \one \Bigl\{ \|(R_{n}+a(R_{n})\rho)\theta + y_i\| \geq R_n \Bigr\} \notag \\
&\times \one \Bigl\{ \bigl( (R_n + a(R_n)\rho)\theta, \, (R_n + a(R_n)\rho)\theta + \by \bigr) \in F_n^{(\ell)} (K,L)\Bigr\}\, h_{t,s}^{(\ell)} (0,\by)\,, \notag
\end{align}
where $S^{d-1}$ is the $(d-1)$-dimensional unit sphere in $\bbr^d$. \\
The following expansion is applied frequently in the following.
For each $i=1,\dots,2k-\ell-1$,
\begin{equation*}
\Bigl|\Bigl|\bigl( R_{n} + a(R_{n})\rho \bigr)\theta + y_i \Bigr|\Bigr| = R_{n} + a(R_{n}) \rho + \langle \theta, y_i \rangle + \gamma_n(\rho, \theta, y_i)\,,
\end{equation*}
so that $\gamma_n(\rho, \theta, y_i) \to 0$ uniformly in $\rho >
0$, $\theta \in S^{d-1}$,  and $\|y_i\| \leq k(t+s)$.

For the application of the dominated convergence theorem, we need
to compute the limit of the expression under the integral sign,
while establishing an integrable upper bound. We first calculate
the limit of the indicator functions. For every $\rho >0$, $\theta
\in S^{d-1}$, and $\| y_i\| \leq k(t+s)$, $i=1,\dots,2k-\ell-1$,
\begin{align*}
\prod_{i=1}^{2k-\ell-1} &\one \Bigl\{ \|(R_{n}+a(R_{n})\rho)\theta + y_i\| \geq R_n \Bigr\}  \\
&\times \one \Bigl\{ \bigl( (R_n + a(R_n)\rho)\theta, \, (R_n + a(R_n)\rho)\theta + \by \bigr) \in F_n^{(\ell)} (K,L)\Bigr\} \\
&\quad \to \one \bigl\{ \by \in E_{K,L}^{(\ell)}(\rho, \theta)\, \bigr\}\,, \ \ \ n\to \infty\,.
\end{align*}

Next, it is clear that for every $\rho >0$, $\bigl( 1 + a(R_n)\rho
/R_n \bigr)^{d-1}$ tends to $1$ as $n \to \infty$ (see
\eqref{e:auxi.slow}) and is bounded above by $2 \bigl( \max \{ 1,
\rho \} \bigr)^{d-1}$.

As for the ratio of the densities in the second line of
\eqref{e:Jk}, we use the basic fact that $1/a$ is flat for $a$,
that is, as $n \to \infty$,
\begin{equation}  \label{e:unif.conv.a}
\frac{a(R_{n})}{a\bigl( R_{n} + a(R_{n}) v \bigr)} \to 1\,, \ \ \text{uniformly on bounded } v \text{-sets};
\end{equation}
see p142 in \cite{embrechts:kluppelberg:mikosch:1997} for details.
Noting that $L$ is also flat for $a$, we have for every $\rho >0$,
\begin{align*}
\frac{f\Bigl( \bigl( R_{n} + a(R_{n}) \rho \bigr) e_1 \Bigr)}{f(R_ne_1)} &= \frac{L\bigl( R_{n} + a(R_{n}) \rho \bigr)}{L(R_n)}\, \exp \Bigl\{ -\psi \bigl( R_{n} + a(R_{n})\rho \bigr) + \psi(R_{n}) \Bigr\} \\
&=\frac{L\bigl( R_{n} + a(R_{n}) \rho \bigr)}{L(R_n)} \exp \Bigl\{ - \int_0^\rho \frac{a(R_n)}{a \bigl( R_n + a(R_n)r \bigr)} dr \Bigr\}  \\
&\to e^{-\rho}, \ \ \text{as } n \to \infty\,.
\end{align*}
To provide an upper bound for the ratio of the densities, let
$\bigl( q_m(n), \, m \geq 0, n \geq 1 \bigr)$ be a sequence
defined by
$$
q_m(n) = a(R_{n})^{-1} \Bigl( \psiinv \bigl(\psi(R_{n}) + m \bigr) - R_{n} \Bigr)\,,
$$
equivalently,
$$
\psi \bigl( R_n + a(R_n) q_m(n) \bigr) = \psi(R_n) + m\,.
$$
Then, for $\epsilon \in \bigl( 0, (d+\gamma (2k-\ell))^{-1}
\bigr)$, there exists an integer $N_{\epsilon} \geq 1$ such that
$$
q_m(n) \leq e^{m\epsilon} / \epsilon \ \ \text{for all } n \geq N_{\epsilon}, m \geq 0\,.
$$
For the proof of this assertion, the reader may refer to Lemma 5.2
in \cite{balkema:embrechts:2004}; see also Lemma 4.7 of
\cite{owada:adler:2015}. Because of the fact that $\psi$ is
non-decreasing, we have, for sufficiently large $n$,
\begin{align*}
&\exp \Bigl\{ -\psi \bigl( R_{n} + a(R_{n})\rho \bigr) + \psi(R_{n}) \Bigr\} \, \one \{\rho > 0 \}  \\
&= \sum_{m=0}^{\infty} \, \one \bigl\{ q_m(n) < \rho \leq q_{m+1}(n) \bigr\}\, \exp \Bigl\{ -\psi \bigl( R_{n} + a(R_{n})\rho \bigr) + \psi(R_{n}) \Bigr\}  \\
&\leq \sum_{m=0}^{\infty} \, \one \bigl\{ 0 < \rho \leq \epsilon^{-1} e^{(m+1)\epsilon} \bigr\}\, e^{-m}.
\end{align*}
Using the bound in \eqref{e:poly.upper},
\begin{align*}
L(&R_{n})^{-1} L \bigl( R_{n} + a(R_{n}) \rho \bigr) \one \{ \rho >0 \}  \leq C \left( 1 + \frac{a(R_{n})}{R_{n}} \rho \right)^{\gamma} \leq
2C \bigl( \max \{ \rho, 1\} \bigr)^{\gamma}.
\end{align*}
Combining these bounds,
$$
\frac{f\Bigl( \bigl( R_{n} + a(R_{n}) \rho \bigr) e_1 \Bigr)}{f(R_ne_1)}\, \one \{ \rho >0 \} \leq 2C \bigl( \max \{ \rho, 1\} \bigr)^{\gamma} \sum_{m=0}^{\infty} \, \one \bigl\{ 0 < \rho \leq \epsilon^{-1} e^{(m+1)\epsilon} \bigr\}\, e^{-m}.
$$
Finally, we turn to
\begin{align*}
\prod_{i=1}^{2k-\ell-1} \frac{f\Bigl( \|(R_{n}+a(R_{n})\rho)\theta + y_i\|e_1 \Bigr)}{f(R_ne_1)} = & \prod_{i=1}^{2k-\ell-1} \frac{L \Bigl( R_{n} + a(R_{n}) \bigl( \rho +  \xi_n(\rho,\theta,y_i) \bigr) \Bigr)}{L(R_n)} \\
&\times \exp \left\{ -\int_0^{\rho + \xi_n(\rho,\theta,y_i)} \frac{a(R_{n})}{a \bigl( R_{n} + a(R_{n} )r \bigr)}\, dr \right\},
\end{align*}
where
$$
\xi_n(\rho, \theta, y) = \frac{\langle \theta, y \rangle + \gamma_n(\rho, \theta, y)}{a(R_{n})}\,.
$$
Since $c=\lim_{n \to \infty}a(R_n)>0$,
$$
A := \sup_{\substack{n \geq 1, \ \rho >0, \\[2pt] \theta \in S^{d-1}, \ \|y\| \leq k(t+s)}} \bigl| \xi_n(\rho, \theta, y) \bigr| < \infty\,.
$$
Therefore, because of the uniform convergence in
\eqref{e:unif.conv.a}, for every $\rho>0$, $\theta \in S^{d-1}$,
and $\| y_i \| \leq k(t+s)$,
$$
\prod_{i=1}^{2k-\ell-1} \frac{f\Bigl( \|(R_{n}+a(R_{n})\rho)\theta + y_i\|e_1 \Bigr)}{f(R_ne_1)} \to \exp \bigl\{ -(2k-\ell-1)\rho - c^{-1} \sum_{i=1}^{2k-\ell-1}\langle \theta, y_i \rangle \bigr\}.
$$
Subsequently, on the set
\begin{align*}
\Bigl\{ &\|(R_{n}+a(R_{n})\rho)\theta + y_i\| \geq R_n\,, \ i=1,\dots,2k-\ell-1 \Bigr\} \\
&=\bigl\{  \rho + \xi_n(\rho, \theta, y_i) \geq 0, \ i=1,\dots,2k-\ell-1 \bigr\}\,,
\end{align*}
we have an obvious upper bound
$$
\prod_{i=1}^{2k-\ell-1} \exp \left\{ -\int_0^{\rho + \xi_n(\rho,\theta,y_i)} \frac{a(R_{n})}{a \bigl( R_{n} + a(R_{n} )r \bigr)}\, dr \right\} \leq 1
$$
from which, together with \eqref{e:poly.upper}, we see that
\begin{align*}
\prod_{i=1}^{2k-\ell-1} \frac{f\Bigl( \|(R_{n}+a(R_{n})\rho)\theta + y_i\|e_1 \Bigr)}{f(R_{n}e_1)} &\leq \prod_{i=1}^{2k-\ell-1} C \left( 1 + \frac{a(R_{n})}{R_{n}} \bigl( \rho + \xi_n(\rho,\theta,y_i) \bigr) \right)^{\gamma} \\
&\leq C^* \bigl( \max \{ \rho,1 \} \bigr)^{\gamma (2k-\ell-1)}.
\end{align*}

From the argument thus far, for every $\rho >0$, $\theta \in
S^{d-1}$, and $\|  y_i \| \leq k(t+s)$, $i=1,\dots,2k-\ell-1$, the
expression under the integral sign in \eqref{e:Jk} eventually
converges to
$$
e^{ -(2k-\ell)\rho - c^{-1} \sum_{i=1}^{2k-\ell-1} \langle \theta,y_i \rangle }\, \one \bigl\{  \by \in E_{K,L}^{(\ell)} (\rho,\theta) \bigr\}\,  h_{t,s}^{(\ell)}(0,\by)\,,
$$
while it possesses an upper bound of the form
$$
C^* \bigl( \max \{ \rho,1 \} \bigr)^{d-1+\gamma (2k-\ell)} \sum_{m=0}^{\infty} \, \one \bigl\{ 0 < \rho \leq \epsilon^{-1} e^{(m+1)\epsilon} \bigr\}\, e^{-m} h_{t,s}^{(\ell)}(0,\by)
$$
for sufficiently large $n$. Because of the restriction in
$\epsilon$, it is elementary to check that
$$
\int_{0}^{\infty}\bigl( \max \{ \rho,1 \} \bigr)^{d-1+\gamma (2k-\ell)} \sum_{m=0}^{\infty} \one \bigl\{ 0 < \rho \leq \epsilon^{-1} e^{(m+1)\epsilon} \bigr\} e^{-m} d \rho  < \infty\,.
$$
As a result of the dominated convergence theorem, we have obtained, as $n\to\infty$,
\begin{align*}
J_k &\sim a(R_n)R_n^{d-1} f(R_ne_1)^{2k-\ell} \int_{S_{d-1}} \hspace{-7pt} J(\theta) d\theta \int_0^\infty d\rho \int_{(\bbr^d)^{2k-\ell-1}} \hspace{-10pt}d\by \\
&\qquad \qquad \times e^{ -(2k-\ell)\rho - c^{-1} \sum_{i=1}^{2k-\ell-1} \langle \theta,y_i \rangle }\, \one \bigl\{  \by \in E_{K,L}^{(\ell)} (\rho,\theta) \bigr\}\,  h_{t,s}^{(\ell)}(0,\by) \\
&=  a(R_n)R_n^{d-1} f(R_ne_1)^{2k-\ell} \ell ! \bigl( (k-\ell)! \bigr)^2 M_{\ell, K, L}(t,s)\,, 
\end{align*}
where the last step follows from the rotation invariance of $h_\cdot$. Hence, we have
\begin{align*}
\text{Cov} &\bigl( G_{n,K,L}(t)\,, G_{n,K,L}(s) \bigr)  \sim \sum_{\ell=1}^k n^{2k-\ell} a(R_n)R_n^{d-1} f(R_ne_1)^{2k-\ell} M_{\ell, K, L}(t,s)\,, \ \ \ n \to \infty\,. \notag
\end{align*}
If $nf(R_ne_1) \to 0$, then the $k$th term in the sum is
asymptotically dominant, and therefore, statement $(i)$ of the
theorem is complete. However, the first term becomes dominant when
$nf(R_ne_1) \to \infty$, in which case, statement $(iii)$ is
established. In addition, if $nf(R_ne_1) \to \xi \in (0,\infty)$,
all the terms in the sum grow at the same rate, and this completes
statement $(ii)$.
\end{proof}
Subsequently, we show the results on finite-dimensional weak
convergence of $\BX_n$ and $(\BX_n^{+}, \BX_n^-)$ defined in
\eqref{e:def.Xn} and \eqref{e:def.Xnpm}, which somewhat parallel
those of Proposition \ref{p:fidi.heavy}. The reader may return to
Section \ref{s:limit.proc.light} to recall the definition and
properties of the limit $(\BW_\ell^+, \BW_\ell^-)$. We omit their
proofs, since the argument in Proposition \ref{p:fidi.heavy} does
apply again with minor modifications.
\begin{proposition}   \label{p:fidi.light}
Assume the conditions of Theorem \ref{t:main.light}. Then, weak
convergences $(i) - (iii)$ in the theorem hold in a
finite-dimensional sense.
Furthermore,  the following results also hold in a finite-dimensional sense. \\
$(i)$ If $nf(R_ne_1) \to 0$ as $n \to \infty$, then
$$
(\BX_n^+, \BX_n^-) \Rightarrow (\BW_k^+, \BW_k^-)\,.
$$
$(ii)$ If $nf(R_ne_1) \to \xi \in (0,\infty)$ as $n \to \infty$, then
$$
(\BX_n^+, \BX_n^-) \Rightarrow \left(\sum_{\ell=1}^k \xi^{2k-\ell} \BW_\ell^+, \, \sum_{\ell=1}^k \xi^{2k-\ell} \BW_\ell^- \right)\,.
$$
$(iii)$ If $nf(R_ne_1) \to \infty$ as $n \to \infty$, then
$$
(\BX_n^+, \BX_n^-) \Rightarrow (\BW_1^+, \BW_1^-)\,.
$$
\end{proposition}
For the same reason as discussed in the preceding subsection, the
next proposition can complete the proof of Theorem
\ref{t:main.light}.
\begin{proposition}
The sequences $(\BX_n^+)$ and $(\BX_n^-)$ are both tight in
$\mathcal D [0,\infty)$, regardless of the limit of $nf(R_ne_1)$.
\end{proposition}
\begin{proof}
We only prove the tightness of $(\BX_n^+)$ but suppress the superscript ``+" from the functions and objects involved during the proof. 
Proceeding completely in the same manner as Proposition
\ref{p:tightness.heavy}, we  have only to show that there exists
$B>0$ such that
\begin{multline}
\tau_n^{-2} \E \Bigl\{ \sum_{\Y_1 \subset \Pn}\sum_{\Y_2 \subset \Pn} \sum_{\Y_3 \subset \Pn}\sum_{\Y_4 \subset \Pn} h_{n,t,s}(\Y_1)\, h_{n,t,s}(\Y_2)\,h_{n,s,r}(\Y_3)\,h_{n,s,r}(\Y_4)\, \label{e:pattern5} \\
\times \one \bigl\{ \text{each } \Y_i \text{ has at least one common elements} \\
\text{with at least one of the other three} \bigr\} \Bigr\} \leq B(t-r)^2
\end{multline}
for all $0 \leq r \leq s \leq t \leq L$ and $n \geq 1$. There are three possibilities to be discussed. \\
$[\text{I}]$ $\ell := |\Y_1 \cap \Y_2| \in \{ 1,\dots,k \}$, $\ellp := |\Y_3 \cap \Y_4| \in \{ 1,\dots,k \}$, and $(\Y_1 \cup \Y_2) \cap (\Y_3 \cup \Y_4) = \emptyset$. \\
$[\text{I} \hspace{-1pt }\text{I}]$ $\ell := |\Y_2 \cap \Y_3| \in \{ 1,\dots,k \}$, $\ellp := |\Y_1 \cap \Y_4| \in \{ 1,\dots,k \}$, and $(\Y_2 \cup \Y_3) \cap (\Y_1 \cup \Y_4) = \emptyset$. \\
$[\text{I}\hspace{-1pt}\text{I}\hspace{-1pt}\text{I}]$. Each
$\Y_i$ has at least one common element with at least one of the
other three, but neither $[\text{I}]$ or $[\text{I}\hspace{-1pt
}\text{I}]$ is true.

Let $B$ be the left hand side of \eqref{e:pattern5}. As for case
$[\text{I}]$, by mimicking the argument in Proposition
\ref{p:tightness.heavy}, we obtain
\begin{align*}
B \leq C^* \tau_n^{-2} n^{4k-\ell-\ellp} a(R_n)^2 R_n^{2(d-1)} f(R_ne_1)^{4k-\ell-\ellp} (t-r)^2 \leq C^* (t-r)^2,
\end{align*}
which proves \eqref{e:pattern5}. Since we can deal with $[\text{I}
\hspace{-1pt }\text{I}]$ in an analogous way, we can turn to case
$[\text{I}\hspace{-1pt}\text{I}\hspace{-1pt}\text{I}]$. Letting
$\ell := 4k-|\Y_1 \cup \Y_2 \cup \Y_3 \cup \Y_4| \in \{ 3,\dots,3k
\}$, the same argument as Proposition \ref{p:tightness.heavy}
yields
\begin{align*}
B \leq C^* \tau_n^{-2} n^{4k-\ell} a(R_n) R_n^{d-1} f(R_ne_1)^{4k-\ell} (t-r)^2 \leq C^* (t-r)^2
\end{align*}
which verifies \eqref{e:pattern5}.
\end{proof}
%\vspace{10pt}

\section{Appendix}

We collect supplemental but important results for the completion
of the main theorems. This result is known as the Palm theory
of Poisson point processes, which is applied a number of times
throughout the proof. 
%The second result is referred to as Stein's
%method for normal approximation, regarding sums of a certain
%weakly dependent variable.
\begin{lemma} (Palm theory for Poisson point processes, \cite{arratia:goldstein:gordon:1989},
Corollary B.2 in \cite{bobrowski:adler:2014}, see also Theorem 1.6 in \cite{penrose:2003})  \label{l:palm1}
Let $(X_i)$ be $\text{i.i.d.}$ $\bbr^d$-valued random variables with common
density $f$. Let $\Pn$ be a Poisson point process on $\bbr^d$ with
intensity $nf$. Let $h(\Y)$, $h_i(\Y)$, $i=1,2,3,4$ be measurable
bounded functions defined for $\Y \in (\bbr^d)^k$. Then,
\begin{align*}
\E \Bigl\{ \sum_{\Y \subset \Pn}  h(\Y) \Bigr\} &= \frac{n^k}{k!} \E \bigl\{ h(\Y) \bigr\}\,,
\end{align*}
and for every $\ell \in \{ 0,\dots,k \}$,
\begin{align*}
\E \Bigl\{ \sum_{\Y_1 \subset \Pn} \sum_{\Y_2 \subset \Pn} h_1(\Y_1)\, h_2(\Y_2)\, \one \bigl\{ |\Y_1 \cap \Y_2| = \ell \bigr\} \Bigr\} &= \frac{n^{2k-\ell}}{\ell ! \bigl( (k-\ell)! \bigr)^2}\, \E \Bigl\{ h_1(\Y_1)\, h_2(\Y_2)\, \one \bigl\{ |\Y_1 \cap \Y_2| = \ell \bigr\} \Bigr\}\,.
\end{align*}
Moreover, for every $\ell_{1}, \ell_{2}, \ell_{3} \in \{ 0,\dots,k
\}$ and $\ell \in \{ 0,\dots,2k \}$, there exists a constant
$C>0$, which depends only on $\ell_{i}$, $\ell$, and $k$ such that
\begin{align*}
&\E \biggl\{  \sum_{\Y_1 \subset \Pn} \sum_{\Y_2 \subset \Pn} \sum_{\Y_3 \subset \Pn} h_1(\Y_1)\, h_2(\Y_2)\, h_3(\Y_3)\, \\
&\quad \quad \quad \quad \times \one \bigl\{ |\Y_1 \cap \Y_2| = \ell_{1}, \, |\Y_1 \cap \Y_3| = \ell_{2}, \, |\Y_2 \cap \Y_3 |= \ell_{3}, \, |\Y_1 \cup \Y_2 \cup \Y_3| = 3k-\ell \bigr\} \biggr\}  \\
&\quad = C n^{3k-\ell} \E \Bigl\{  h_1(\Y_1)\, h_2(\Y_2)\, h_3(\Y_3)\,\\
&\quad \quad \quad \quad \times \one \bigl\{ |\Y_1 \cap \Y_2| = \ell_{1}, \, |\Y_1 \cap \Y_3| = \ell_{2}, \, |\Y_2 \cap \Y_3| = \ell_{3}, \, |\Y_1 \cup \Y_2 \cup \Y_3| = 3k-\ell \bigr\} \Bigr\}\,.
\end{align*}
Similarly, for $\ell_{i,j} \in \{ 0,\dots,k \}$, $m_{p,q,r} \in \{
0,\dots,k \}$, $i,j,p,q,r \in \{1,2,3,4\}$ with $i \neq j, p \neq
q, p \neq r, q \neq r$, and $\ell \in \{0,\dots,3k\}$, there
exists a constant $C>0$, which depends only on $\ell_{i,j}$,
$m_{p,q,r}$, $\ell$, and $k$ such that
\begin{align*}
&\E \biggl\{  \sum_{\Y_1 \subset \Pn} \sum_{\Y_2 \subset \Pn} \sum_{\Y_3 \subset \Pn} \sum_{\Y_4 \subset \Pn} h_1(\Y_1)\, h_2(\Y_2)\, h_3(\Y_3)\, h_4(\Y_4)\\
&\quad \quad \quad \times \one \bigl\{ |\Y_i \cap \Y_j| = \ell_{i,j}, \, i,  j \in  \{ 1,2,3,4 \}, i \neq j\,, \\
&\quad \quad \quad \quad \quad  \quad |\Y_p \cap \Y_q  \cap \Y_r| = m_{p,q,r}, \, p,q,r \in \{ 1,2,3,4 \}, \, p \neq q, p \neq q, q \neq r, \\
&\quad \quad \quad \quad \quad \quad \quad |\Y_1 \cup \Y_2 \cup \Y_3 \cup \Y_4| = 4k-\ell \bigr\} \biggr\}  \\
&\quad = C n^{4k-\ell} \E \Bigl\{  h_1(\Y_1)\, h_2(\Y_2)\, h_3(\Y_3)\, h_4(\Y_4)\\
&\quad \quad \quad \quad \quad \times \one \bigl\{ |\Y_i \cap \Y_j| = \ell_{i,j}, \, i,  j \in  \{ 1,2,3,4 \}, i \neq j\,, \\
&\quad \quad \quad \quad \quad \quad \quad  \quad |\Y_p \cap \Y_q  \cap \Y_r| = m_{p,q,r}, \, p,q,r \in \{ 1,2,3,4 \}, \, p \neq q, p \neq q, q \neq r, \\
&\quad \quad \quad \quad \quad \quad \quad \quad \quad |\Y_1 \cup \Y_2 \cup \Y_3 \cup \Y_4| = 4k-\ell \bigr\} \biggr\}\,.
\end{align*}
\end{lemma}

%Before presenting Stein's method for normal approximation, we
%recall the definition of a \textit{dependency graph}. Given a
%graph $(I,E)$, we denote $i \sim j$ for $i,j \in I$ iff there is
%an edge connecting $i$ and $j$. Let $(Y_i, \, i \in I)$ be a
%family of random variables. Then, $(I,\sim)$ is said to be a
%dependency graph with respect to $(Y_i, \, i \in I)$ if for any
%two disjoint subsets $I_1, I_2$ of $I$ with no edges connecting
%$I_1$ and $I_2$, $(Y_i, \, i \in I_1)$ is independent of $(Y_i, \,
%i \in I_2)$.
%\begin{lemma}(Stein's normal approximation, Theorem 2.4 in \cite{penrose:2003}) \label{l:stein.normal}
%Let $(Y_i, \, i \in I)$ be a finite collection of random variables
%with $\E \{Y_i \} = 0$ for each $i \in I$ and $\E \Bigl\{ \bigl(
%\sum_{i \in I} Y_i \bigr)^2 \Bigr\}=1$. Let $(I,\sim)$ be a
%dependency graph with respect to $(Y_i, \, i \in I)$ with maximum
%degree $D-1$. Let $\Phi$ be the distribution function of a
%standard normal random variable. Then, for all $\lambda>0$,
%\begin{align*}
%\biggl| \P \Bigl\{ \sum_{i \in I} Y_i \leq \lambda \biggr\} - \Phi(\lambda) \Bigr| &\leq 2(2\pi)^{-1/4} \sqrt{D^2 \sum_{i \in I} \E \bigl\{ |Y_i|^3 \bigr\}}  +6 \sqrt{D^3 \sum_{i \in I} \E %\bigl\{ |Y_i|^4 \bigr\}}.
%\end{align*}
%\end{lemma}

%{\bf Acknowledgment} The author should like to express his
%gratitude to Professor Robert Adler of Technion-Israel Institute
%of Technology for fruitful discussions and his sets of helpful
%comments throughout this research.

%\bibliography{Takashi_ref}

\begin{thebibliography}{10}

\bibitem{adler:bobrowski:weinberger:2014}
R.~J. Adler, O.~Bobrowski, and S.~Weinberger.
\newblock Crackle: The homology of noise.
\newblock {\em Discrete \& Computational Geometry}, 52:680--704, 2014.

\bibitem{arratia:goldstein:gordon:1989}
R.~Arratia, L.~Goldstein, and L.~Gordon.
\newblock Two moments suffice for poisson approximations: the chen-stein
  method.
\newblock {\em The Annals of Probability}, 17:9--25, 1989.

\bibitem{balkema:embrechts:2004}
G.~Balkema and P.~Embrechts.
\newblock Multivariate excess distributions.
\newblock {www.math.ethz.ch/~embrecht/ftp/guuspe08Jun04.pdf}, 2004.

\bibitem{balkema:embrechts:2007}
G.~Balkema and P.~Embrechts.
\newblock {\em High Risk Scenarios and Extremes: A Geometric Approach}.
\newblock European Mathematical Society, 2007.

\bibitem{balkema:embrechts:nolde:2010}
G.~Balkema, P.~Embrechts, and N.~Nolde.
\newblock Meta densities and the shape of their sample clouds.
\newblock {\em Journal of Multivariate Analysis}, 101:1738--1754, 2010.

\bibitem{balkema:embrechts:nolde:2013}
G.~Balkema, P.~Embrechts, and N.~Nolde.
\newblock The shape of asymptotic dependence.
\newblock {\em Springer Proceedings in Mathematics \& Statistics, Special
  volume "Prokhorov and Contemporary Probability Theory"}, 33:43--67, 2013.

\bibitem{bhattacharya:ghosh:1992}
R.~N. Bhattacharya and J.~K. Ghosh.
\newblock A class of $u$-statistics and asymptotic normality of the number of
  $k$-clusters.
\newblock {\em Journal of Multivariate Analysis}, 43:300--330, 1992.

\bibitem{billingsley:1999}
P.~Billingsley.
\newblock {\em Convergence of Probability Measures, \rm{2nd edition}}.
\newblock Wiley, New York, 1999.

\bibitem{bobrowski:adler:2014}
O.~Bobrowski and R.~J. Adler.
\newblock Distance functions, critical points, and topology for some random
  complexes.
\newblock {\em Homology, Homotopy and Applications}, 16:311--344, 2014.

\bibitem{bobrowski:kahle:2014}
O.~Bobrowski and M.~Kahle.
\newblock Topology of random geometric complexes: a survey.
\newblock {arXiv:1409.4734}, 2014.

\bibitem{chen:jia:2001}
X.~Chen and X.~Jia.
\newblock Package routing algorithms in mobile ad-hoc wireless networks.
\newblock {\em Proceeding of the Workshop on Wireless Networks and Mobile
  Computing in conjunction with the 2001 International Conference on Parallel
  Processing}, pages 485--490, 2001.

\bibitem{dabrowski:dehling:mikosch:sharipov:2002}
A.~R. Dabrowski, H.~G. Dehling, T.~Mikosch, and O.~Sharipov.
\newblock Poisson limits for $u$-statistics.
\newblock {\em Stochastic Processes and their Applications}, 99:137--157, 2002.

\bibitem{dehaan:ferreira:2006}
L.~{\rm de Haan} and A.~Ferreira.
\newblock {\em Extreme Value Theory: An Introduction}.
\newblock Springer, New York, 2006.

\bibitem{embrechts:kluppelberg:mikosch:1997}
P.~Embrechts, C.~Kl\"uppelberg, and T.~Mikosch.
\newblock {\em Modelling Extremal Events: for Insurance and Finance}.
\newblock Springer, New York, 1997.

\bibitem{gilbert:1961}
E.~N. Gilbert.
\newblock Random plane networks.
\newblock {\em Journal of the Society for Industrial and Applied Mathematics},
  9:533--543, 1961.

\bibitem{hafner:1972}
R.~Hafner.
\newblock The asymptotic distribution of random clumps.
\newblock {\em Computing}, 10:335--351, 1972.

\bibitem{hekmat:2006}
R.~Hekmat.
\newblock {\em Ad-hoc Networks: Fundamental Properties and Network Topologies}.
\newblock Springer, New York, 2006.

\bibitem{kahle:2011}
M.~Kahle.
\newblock Random geometric complexes.
\newblock {\em Discrete \& Computational Geometry}, 45:553--573, 2011.

\bibitem{kahle:meckes:2013}
M.~Kahle and E.~Meckes.
\newblock Limit theorems for betti numbers of random simplicial complexes.
\newblock {\em Homology, Homotopy and Applications}, 15:343--374, 2013.

\bibitem{owada:adler:2015}
T.~Owada and R.~J. Adler.
\newblock Limit theorems for point processes under geometric constraints (and
  topological crackle).
\newblock {arXiv:1503.08416}, 2015.

\bibitem{penrose:2003}
M.~Penrose.
\newblock {\em Random Geometric Graphs, Oxford Studies in Probability 5}.
\newblock Oxford University Press, Oxford, 2003.

\bibitem{resnick:1987}
S.~Resnick.
\newblock {\em Extreme Values, Regular Variation and Point Processes}.
\newblock Springer-Verlag, New York, 1987.

\bibitem{resnick:2007}
S.~Resnick.
\newblock {\em Heavy-Tail Phenomena: Probabilistic and Statistical Modeling}.
\newblock Springer, New York, 2007.

\bibitem{schulte:thale:2012}
M.~Schulte and C.~Th\"{a}le.
\newblock The scaling limit of poisson-driven order statistics with
  applications in geometric probability.
\newblock {\em Stochastic Processes and their Applications}, 122:4096--4120,
  2012.

\bibitem{silverman:brown:1978}
B.~Silverman and T.~Brown.
\newblock Short distances, flat triangles and poisson limits.
\newblock {\em Journal of Applied Probability}, 15:815--825, 1978.

\bibitem{stojmenovic:seddigh:zunic:2002}
I.~Stojmenovic, M.~Seddigh, and J.~Zunic.
\newblock Dominating sets and neighbor elimination-based broadcasting
  algorithms in wireless networks.
\newblock {\em IEEE Transactions on Parallel and Distributed Systems},
  13:14--25, 2002.

\bibitem{weber:1983}
N.~C. Weber.
\newblock Central limit theorems for a class of symmetric statistics.
\newblock {\em Mathematical Proceedings of the Cambridge Philosophical
  Society}, 94:307--313, 1983.

\bibitem{yogeshwaran:subag:adler:2014}
D.~Yogeshwaran, E.~Subag, and R.~J. Adler.
\newblock Random geometric complexes in the thermodynamic regime.
\newblock {\em Probability Theory and Related Fields}, 2016.
\newblock {In press, arXiv:1403.1164}.

\end{thebibliography}

\end{document}